\documentclass[review,onefignum,onetabnum,pdftex]{siamonline190516}

\usepackage{lipsum}
\usepackage{amsfonts}
\usepackage{graphicx}
\usepackage{epstopdf}
\usepackage{algorithmic}
\ifpdf
  \DeclareGraphicsExtensions{.eps,.pdf,.png,.jpg}
\else
  \DeclareGraphicsExtensions{.eps}
\fi

\usepackage{enumitem}
\setlist[enumerate]{leftmargin=.5in}
\setlist[itemize]{leftmargin=.5in}

\newsiamremark{remark}{Remark}
\newsiamremark{hypothesis}{Hypothesis}
\crefname{hypothesis}{Hypothesis}{Hypotheses}
\newsiamthm{claim}{Claim}

\headers{Augmented Lagrangian Sketching}{Md Sarowar Morshed}

\title{ALS: Augmented Lagrangian Sketching Methods for Linear Systems}

\author{Md Sarowar Morshed \thanks{Independent Researcher
  (\email{riponsarowar@outlook.com}, \url{https://sarowar.netlify.app/}).}}

\usepackage{amsopn}

\usepackage{wrapfig}
\usepackage{hyperref}

\usepackage{subcaption}
\usepackage{float}

\usepackage{cite}

\usepackage{amsfonts}
\usepackage{amssymb,amsmath,latexsym,tensor}
\usepackage{mathrsfs}
\usepackage{savesym}
\usepackage{graphicx}
\usepackage{adjustbox}
\usepackage[pdftex]{color}
\usepackage[font=footnotesize,labelfont=bf]{caption}
\savesymbol{AND}
\usepackage{epstopdf}
\usepackage{algorithm}
\usepackage{algorithmic}

\newcommand{\R}{\mathbb{R}}

\DeclareMathOperator{\E}{\mathbb{E}}

\DeclareMathOperator*{\argmin}{arg\,min}
\DeclareMathOperator*{\argmax}{arg\,max}

\usepackage{mathtools}
\DeclarePairedDelimiterX{\norm}[1]{\lVert}{\rVert}{#1}
\DeclarePairedDelimiterX{\inp}[2]{\langle}{\rangle}{#1, #2}

\usepackage{empheq}

\newtheorem{assumption}{Assumption}

\usepackage[makeroom]{cancel}
\usepackage{mdframed}

\usepackage{array,multirow}

\usepackage{amssymb}

\usepackage[utf8]{inputenc}

\usepackage{enumitem}

\ifpdf
\hypersetup{
  pdftitle={ALS: Augmented Lagrangian Sketching Methods for Linear Systems},
  pdfauthor={Md Sarowar Morshed}
}
\fi

\begin{document}

\maketitle

\begin{abstract}
We develop two fundamental stochastic sketching techniques \textemdash \\ \textit{Penalty Sketching} (PS) and \textit{Augmented Lagrangian Sketching} (ALS) for solving consistent linear systems. The proposed PS and ALS techniques extend and generalize the scope of \textit{Sketch \& Project} (SP) method by introducing Lagrangian penalty sketches. In doing so, we recover SP methods as special cases and furthermore develop a family of new stochastic iterative methods. By varying sketch parameters in the proposed PS method, we recover novel stochastic methods such as \textit{Penalty Newton Descent}, \textit{Penalty Kaczmarz}, \textit{Penalty Stochastic Descent}, \textit{Penalty Coordinate Descent}, \textit{Penalty Gaussian Pursuit}, and \textit{Penalty Block Kaczmarz}. Furthermore, the proposed ALS method synthesizes a wide variety of new stochastic methods such as \textit{Augmented Newton Descent}, \textit{Augmented Kaczmarz}, \textit{Augmented Stochastic Descent}, \textit{Augmented Coordinate Descent}, \textit{Augmented Gaussian Pursuit}, and \textit{Augmented Block Kaczmarz} into one framework. Moreover, we show that the developed PS and ALS frameworks can be used to reformulate the original linear system into equivalent stochastic optimization problems namely \textemdash \textit{Penalty Stochastic Reformulation} and \textit{Augmented Stochastic Reformulation}. We prove global convergence rates for the PS and ALS methods as well as sub-linear $\mathcal{O}(\frac{1}{k})$ rates for the Cesaro average of iterates. The proposed convergence results hold for a wide family of distributions of random matrices, which provides the opportunity of fine-tuning the randomness of the method suitable for specific applications. Finally, we perform computational experiments that demonstrate the efficiency of our methods compared to the existing SP methods.
\end{abstract}



\section{Introduction}
\label{als:ls:sec:intro}
We develop novel iterative sketching method for solving the affine projection problem given by. 
\begin{align}
\label{als:ls:1}
  x^* = \argmin_{x \in \R^{n}} \|x-x_0\|^2_B \quad \textbf{subject to} \quad  Ax = b,
\end{align}
where, $B \in \R^{n \times n}$ is a symmetric positive definite matrix and $ x_0 \in \R^{n}$ is a given vector. We use $\|x\|_B^2 := \langle Bx,x \rangle$ to denote a weighted norm. 
When the system $Ax = b$ has a unique solution, problem \eqref{als:ls:1} reduces to one of the fundamental problems of \textit{Numerical Linear Algebra}. Problem \eqref{als:ls:1} is also central to a wide range of interdisciplinary areas such as \textit{ Convex Optimization}, \textit{Numerical Linear Algebra}, \textit{Scientific Computing}, \textit{Computer Science},  \textit{Machine Learning}, and \textit{Management Science}, etc. In \textit{Machine Learning}, variants of problem \eqref{als:ls:1} found applications in subareas such as \textit{Least-Square Support Vector Machines} \cite{SVM:2007}, \textit{Gaussian Markov Random Fields} \cite{GaussianMR}, \textit{Gaussian Processes} \cite{rasmussen_williams_2008},  graph-based \textit{Semi-Supervised Learning} and \textit{Graph-Laplacian} problems \cite{bengio:2006}, etc. Even though Krylov-type methods are industry standard for solving problem \eqref{als:ls:1} in a low-dimensional setting. However, in a huge dimensional setup, direct methods are infeasible as they require a huge amount of matrix-vector multiplications. Recent works suggest that randomized projection methods are more efficient \cite{strohmer:2008,lewis:2010} owing to their low storage cost and flops complexity.

\subsection{Related Works}
The most simplest and oldest projection methods is the so-called Kaczmarz method discovered in 1937 \cite{kaczmarz:1937}. Originally Kaczmarz proposed the cyclic Kaczmarz method that was deterministic in nature and remained unnoticed to the research community. It was rediscovered in the 1980s in the form of \textit{Algebraic Reconstruction Techniques} (ART) in the area of image reconstruction \cite{gordon:1970}. Subsequently, afterward, it found applications in multidisciplinary areas such as \textit{Computer Tomography} \cite{Censor:1988,herman:2009}, \textit{Signal Processing} \cite{lorenz:2014,quadratic:2016}, \textit{Distributed Computing} \cite{elble:2010,fabio:2012}, etc. Research into Kaczmarz-type methods boomed after the seminal work of Strohmer \textit{et. al} \cite{strohmer:2008}. They showed that instead of a cyclic manner one can design efficient variants of the Kaczmarz algorithms by choosing projection planes randomly. Up until then, several extensions and generalizations of the Kaczmarz method have been proposed for solving a wide range of problems such as \textit{Linear System}, \textit{Linear Feasibility}, \textit{Least-Square}, \textit{Quadratic Equations}, etc \cite{lewis:2010,needell:2010,eldar:2011,zouzias:2013,lee:2013,NEEDELL:2014,agaskar:2014,ma:2015,NEEDELL:2015, blockneddel:2015,needell:2016, quadratic:2016, wright:2016,needell:2016,haddock:2017, greedbai:2018,razaviyayn:2019, needell2019block}.

\paragraph{\textit{Sketch \& Project} Method \& Extensions.} Recently, it has been shown that one can synthesize almost all of the available projection-type iterative methods for solving \eqref{als:ls:1} into one framework namely the SP method \cite{gower:2015}. Gower \textit{et. al} showed that by varying a positive definite matrix $B$ and sketching matrix $S$ in the SP method, one can obtain a family of randomized iterative methods such as \textit{Randomized Newton}, \textit{Randomized Kaczmarz}, \textit{Randomized Coordinate Descent}, \textit{Random Gaussian Pursuit}, and \textit{Randomized Block Kaczmarz}, etc. Subsequently, after that, they extended the scope of the SP method to various well-known methods such as \textit{Quasi-Newton} \cite{gower:2017,NIPS:2018}, \textit{Newton} \cite{gower2019rsn}, \textit{Newton-Raphson} \cite{yuan2020sketched}, etc. For instance, in their work \cite{gower2016linearly} Gower \textit{et. al} showed that one can interpret available \textit{Quasi-Newton} methods such as \textit{Davidon–Fletcher–Powell}, \textit{Powell-Symmetric-Broyden}, \textit{Bad Broyden}, and \textit{Broyden–Fletcher–Goldfarb–Shanno} as special cases of the SP framework. They also provided a stochastic dual ascent framework of the SP method \cite{gower2016stochastic}. Their works prompted several generalizations and acceleration of the SP method for solving LS \cite{loizou:2017, NIPS:2018,richtrik2017stochastic,gower2019adaptive}, \textit{Pseudoinverse} \cite{gower2016linearly},
\textit{Convex Feasibility} \cite{necoara:2019}, \textit{Linear Feasibility} \cite{morshed:sketching}, \textit{Variance Reduction} \cite{sega:2018,Gower2020} problems. One of the key steps in the SP method is the choice of sketching matrix $S$ at each iteration. So far available sampling strategies in the broader context of the SP method can be categorized into the following: \textit{Uniform Sampling} \cite{strohmer:2008, lewis:2010,wright:2016}, \textit{Maximum Distance Sampling} \cite{motzkin,nutini:2016}, \textit{Kaczmarz Motzkin Sampling} \cite{haddock:2017, haddock:2019, Morshed2019, morshed2020generalization,morshed:momentum,morshed:sketching, morshed2020stochastic}, Capped Sampling \cite{bai:2018,gower2019adaptive,morshed2020stochastic}. Recently Gower \textit{et. al} \cite{gower2019adaptive} and Morshed \textit{et. al} \cite{morshed2020stochastic} combined all of the available sampling strategies under some greedy sampling frameworks.

\subsection{Notation}
We follow standard linear algebra notation throughout the paper. The notation  
$\mathbb{R}^{m\times n}$ will be used to denote the set of real-valued matrices. Let $A \in \mathbb{R}^{m\times n} $ be any arbitrary matrix, then by $a_i^\top $, $A_{j}$, $A^{\dagger}$, $\|A\|_F$, $\textbf{Rank}(A)$, $\textbf{Range}(A)$, $\textbf{Null}(A)$, $\lambda_{\min}^+(A)$, $\lambda_{\max}(A)$ we denote the $i^{th}$ row, the $j^{th}$ column, the Moore-Penrose pseudo-inverse, Frobenius norm, rank, range, null, the smallest nonzero eigenvalue, the largest eigenvalue of matrix $A$, respectively. For any positive definite matrix $B$, we denote $\langle x, y \rangle_B = x^\top By $. In the following, we define sketching matrices that will be used throughout the paper.

\begin{definition}
\label{als:ls:d:1}
Let $\mathcal{D}$ be a distribution over real $m \times p$ matrices (discrete/ continuous). We call random matrix $S \sim \mathcal{D}$ as the sketching matrix and $p$ as the sketch size. We will denote the resulting expectation as $\E[\cdot \ | \ S \sim \mathcal{D} ] = \E_{\mathcal{D}}[\ \cdot \ ]$.
\end{definition}

\paragraph{Contributions.}
In this work, we made the following fundamental algorithmic contributions in the theory of sketching methods:
\begin{itemize}
    \item We propose two sketching methods namely the \textit{Penalty Sketching} (PS) and the \textit{Augmented Lagrangian Sketching} (ALS) that are based on the idea of the Lagrangian penalty function.
    \item The proposed sketching techniques generalize the SP method and extend the scope of the SP framework.
    \item Based on the proposed sketching techniques, we provide two 
    novel stochastic reformulations of the LS problem that extends the scope of the stochastic reformulation proposed in \cite{richtrik2017stochastic}.
    \item We design efficient variants of the proposed algorithms by introducing several greedy sketching rules.
\end{itemize}
Moreover, we establish global convergence results of the proposed PS \& ALS algorithms. We also prove convergence for the Cesaro average of iterates, i.e., $\bar{x}_k = 1/k \sum \nolimits_{i =0}^k x_i $ generated by both algorithms. From our convergence results, one can obtain existing results as well as new results by varying parameters such as penalty, sketching matrix, etc. In the following table, we summarize special algorithms and the variants that can be recovered from the proposed algorithms:

\paragraph{Outline.}
In section \ref{als:ls:sec:lp}, we provide technical details of the SP method and introduce the proposed sketching methods. At the end of section \ref{als:ls:sec:lp}, we discuss the properties of the proposed sketching methods as well as their connection to existing methods. In section \ref{als:ls:sec:conv}, we provide the convergence results of the proposed algorithms. Various convenient and greedy sampling techniques are discussed in section \ref{als:ls:sec:sampl}. In section \ref{als:ls:sec:sp}, we derive and discuss various algorithms that can be obtained as special cases from the proposed methods. To validate the efficiency of the proposed methods, we carry out numerical experiments in section \ref{als:ls:sec:num}. We conclude the paper with the conclusion section.

\section{Lagrangian Penalty Approaches}
\label{als:ls:sec:lp}
We first discuss the existing \textit{Sketch \& Project} (SP) framework for solving problem (1). In doing so, we will  also point out some issues with the SP method. We then motivate the new Lagrangian penalty-based approaches for solving problem \eqref{als:ls:1}.

\paragraph{The Sketch-and-Project method.} The ingredients for the SP method include 1) Sketching matrix $S \sim \mathcal{D}$, 2) Positive definite matrix $B$. Given a random iterate $x_k$, the SP method finds the closest point to $x_k$ in such a way that $x_{k+1}$ solves the sketched problem, i.e.,
\begin{align}
\label{als:ls:2}
x_{k+1} =     \argmin_x \frac{1}{2}\|x-x_k\|^2_B \quad \textbf{s.t} \quad  S^\top Ax = S^\top  b.
\end{align}
The closed form solution of the above problem is given by
\begin{align}
\label{als:ls:3}
x_{k+1} = x_k - B^{-1} A^\top S (S^\top AB^{-1}A^\top S)^{\dagger} S^\top (Ax_k-b).
\end{align}
As detailed in~\cite{richtrik2017stochastic}, the SP method can also be viewed as a type of stochastic gradient descent method for solving 
\begin{align}
 \label{als:ls:4}
    x = \argmin_{x} \|Ax-b\|^2_{\E_{\mathcal{D}}[S(S^\top AB^{-1}A^\top S)^{\dagger}S^\top ]}.
\end{align}
In~\cite{richtrik2017stochastic} the authors showed that under mild technical conditions\footnote{For instance if $\E[S(S^\top AB^{-1}A^\top S)^{\dagger}S^\top] \succ 0$ then \eqref{als:ls:1} and \eqref{als:ls:4} are equivalent.} on the sketching matrices the reformulation is equivalent to \eqref{als:ls:1}.

One draw back of the SP method is that it is \emph{discontinuous} with respect to the data $(A,b).$  Indeed, the pseudoinverse in~\eqref{als:ls:3} is only continuous over invertible matrix. When the sketching dimension $p$ is large, the matrix within the pseudo-inverse may not be invertible. Consequently, the SP method may be sensitive to mispecification in the data or even small numerical inaccuracies. The discontinuity is a result of the hard constraint~\eqref{als:ls:3}. To side step this issue, we propose new penalty based sketching methods.

Indeed, one can solve the optimization problem of \eqref{als:ls:2} using several techniques. Our main idea is to
treat problem \eqref{als:ls:2} as a constrained minimization problem and then solve it using Lagrangian penalty-based techniques. The main ingredients for the Lagrangian penalty-based methods include: 1) two positive definite matrices $G \in \R^{p \times p}, B \in \R^{n \times n} $, 2) sketching matrix $S \sim \mathcal{D}$, and 3) penalty parameter $\rho > 0$.

\subsection{Penalty Sketching (PS)} 
Instead of imposing a hard constraint as done in the SP method, the PS (Penalty Sketching) method introduces a penalty term in the objective function of \eqref{als:ls:2} that corresponds to the constraint violation. We  then solve this unconstrained optimization problem to find the next update. Our intuition is that with the above setup, we can force the algorithm to minimize the constraint violation by increasing the penalty gradually. Introducing the constraints of problem \eqref{als:ls:2} as a penalty, we get the following optimization problem
\begin{align}
\label{als:ls:5}
x_{k+1} =     \argmin_x \mathcal{L}(x, \rho) =  \frac{1}{2} \|x-x_k\|^2_B + \frac{\rho}{2}  \big \| S^\top (Ax-b) \big \|^2_{G^{-1}}.
\end{align}

We can explicitly derive the closed-form solution of the above optimization problem. The solution is unique and provided in the next lemma.

\begin{lemma}
\label{als:ls:lem:0}
The following update solves the optimization problem in \eqref{als:ls:5}:
\begin{align*}
x_{k+1}= x_k -  B^{-1}  A^\top  S \left( \frac{1}{\rho} G+ S^\top A B^{-1} A^\top S\right)^{-1} S^\top  (A x_k - b).
\end{align*}
\end{lemma}

\begin{proof}
Setting $\frac{\partial \mathcal{L}(x, \rho)}{\partial x} = 0$, we get
\begin{align*}
    0 = \frac{\partial \mathcal{L}(x, \rho)}{\partial x} = B(x-x_k)+  \rho A^\top  S G^{-1} S^\top (Ax-b).
\end{align*}
Isolating $x$ give the following solution $x_{k+1}$:
\begin{align}
\label{als:ls:6}
x_{k+1}= \left(B+ \rho A^\top SG^{-1}S^\top A \right)^{-1} \left[ Bx_k + \rho A^\top S G^{-1} S^\top b \right].
\end{align}
The solution is unique as the matrix $B+ \rho A^\top SG^{-1}S^\top A $ is invertible for any $\rho > 0$. Since, $B, G \succ 0$, by Woodburry matrix identity we have the following:
\begin{align*}
 \left(B+ \rho A^\top S G^{-1} S^\top A \right)^{-1} = B^{-1} -  B^{-1}  A^\top  S \left( \frac{G}{\rho}+ S^\top A B^{-1} A^\top S\right)^{-1} S^\top  A B^{-1}.
\end{align*}
Using the above, we can simplify the first term of \eqref{als:ls:6} as follows:
\begin{align}
\label{als:ls:7}
  & \left(B+ \rho A^\top SG^{-1}S^\top A \right)^{-1}  Bx_k \nonumber \\
  & \qquad  = x_k - B^{-1}  A^\top  S \left( \frac{1}{\rho} G + S^\top A B^{-1} A^\top S\right)^{-1} S^\top  A x_k.
\end{align}
Similarly the second term of \eqref{als:ls:6} can be simplified as follows:
\begin{align}
\label{als:ls:8}
 & \left(B+ \rho A^\top S G^{-1} S^\top A \right)^{-1} A^\top S G^{-1} \nonumber \\
 & = B^{-1}  A^\top S G^{-1}-  B^{-1}  A^\top  S \left( \frac{1}{\rho} G+ S^\top A B^{-1} A^\top S\right)^{-1} S^\top  A B^{-1}  A^\top SG^{-1} \nonumber \\
 & =  B^{-1}  A^\top SG^{-1} -  B^{-1}  A^\top  S \left[ I - \frac{1}{\rho}\left( \frac{1}{\rho} G+ S^\top A B^{-1} A^\top S\right)^{-1} G\right]G^{-1} \nonumber \\
 & = \frac{1}{\rho} B^{-1}  A^\top S \left( \frac{1}{\rho} G+ S^\top A B^{-1} A^\top S\right)^{-1}.
\end{align}
Substituting \eqref{als:ls:7} and \eqref{als:ls:8} in \eqref{als:ls:6}, we get the following:
\begin{align}
\label{als:ls:9}
x_{k+1}= x_k -  B^{-1}  A^\top  S \left( \frac{1}{\rho} G+ S^\top A B^{-1} A^\top S\right)^{-1} S^\top  (A x_k - b).
\end{align}
This proves the lemma.
\end{proof}

 \begin{algorithm}
\caption{$x_{K+1} = \textbf{PS}(A,b, B, G, c, \mathcal{D}, K)$}
\label{als:ls:alg:ps}
\begin{algorithmic}
\STATE{Choose initial point $x_0 \in \R^n,  \ \rho_0 \in \R$}
\WHILE{$k \leq K$}
\STATE{Draw a fresh sample $S \sim \mathcal{D}$ and update
\begin{align*}
    & x_{k+1} = x_k -  B^{-1}  A^\top  S \left( \frac{1}{\rho_k} G+ S^\top A B^{-1} A^\top S\right)^{-1} S^\top  (A x_k - b) ; \\
    & \rho_{k+1} = c \rho_k ;
\end{align*}
$k \leftarrow k+1$;}
\ENDWHILE
\end{algorithmic}
\end{algorithm}

\begin{remark}
Note that, the optimization problem \eqref{als:ls:5} is equivalent to solving $$\min g(x) = \E[g_S(x)], \mbox{  where  } \ g_S(x) = \frac{1}{2}\norm{Ax-b}^2_{\E[SG^{-1}S^\top ]},$$ using the stochastic proximal point method. Indeed, we can write the update \eqref{als:ls:5} as follows:
\begin{align*}
  x_{k+1} = \argmin_x \left\{g_S(x) + \frac{1}{2 \rho} \norm{x-x_k}_B^2\right\}  = \textbf{prox}^B_{ \rho \ g_{S}(x_k)}
\end{align*}
here, we define the proximal operator of a function $\psi : \R^n \rightarrow \R$ with respect to the $B-$norm as $\textbf{prox}^B_{\psi}(y) = \argmin_{x \in \R^n} \left\{ \psi(x) + \frac{1}{2} \norm{x-y}_B^2\right\}$. Also, note that the above proximal viewpoint is different than the one provided in \cite{richtrik2017stochastic} for the SP method. In \cite{richtrik2017stochastic}, the authors provided the following proximal method viewpoint:
\begin{align*}
  x_{k+1} = \argmin_x \left\{h_S(x) + \frac{1-\omega}{2 \omega} \norm{x-x_k}_B^2\right\}  = \textbf{prox}^B_{\frac{\omega}{ 1- \omega} h_{S}(x_k)}
\end{align*}
where, $h_{S}(x) = \frac{1}{2} \|Ax-b\|^2_{S(S^\top AB^{-1}A^\top S)^{\dagger}S^\top }$. 
\end{remark}

\subsection{Augmented Lagrangian Sketching (ALS)}
We now present the Augmented Lagrangian Sketching (ALS) framework for solving \eqref{als:ls:1}. The ALS framework introduces the dual variable $z$ that works as an estimator to the Lagrange multiplier. Following the Lagrangian dual formulation, the optimization problem of \eqref{als:ls:2} can be reformulated as follows:
\begin{align}
\label{als:ls:001}
 \min_x \max_z  \mathcal{L}(x,z) \overset{\textbf{def}}{:=}  \frac{1}{2} \|x-x_k\|^2_B + z^\top  S^\top (Ax-b).
\end{align}
Since the function $\mathcal{L}(x,z)$ is non-smooth with respect to the dual variable $z$, we regularize the objective using the distance of $z$ to a given previous iterate $z_k$ as follows:
\begin{align}
\label{als:ls:002}
 (x_{k+1},z_{k+1}) = \argmin_x \argmax_z  \ \mathcal{L}(x,z) - \frac{1}{2 \rho} \underbrace{\|z-z_k\|^2_G}_{\textbf{smoothing term}}.
\end{align}
This now makes the objective smooth with respect to $z$, and also allows us to uniquely define the next dual iterate given by
\begin{align}
\label{als:ls:003}
    z_{k+1} & = \argmax_z \ \mathcal{L}(x,z)  - \frac{1}{2 \rho} \|z-z_k\|^2_G = z_k + \rho \ G^{-1} S^\top  (Ax-b).
\end{align} 
Now, substituting $z = z_{k+1} = z_k + \rho \ G^{-1} S^\top  (Ax-b)$ in the objective function of \eqref{als:ls:002}, we estimate the primal update as follows:
\begin{align}
\label{als:ls:004}
    x_{k+1} & = \argmin_x \ \mathcal{L}(x,z_{k+1})  - \frac{1}{2 \rho} \|z_{k+1}-z_k\|^2_G \nonumber \\
    & = \argmin_x \ \frac{1}{2} \|x-x_k\|^2_B + z_{k+1}^\top  S^\top (Ax-b) - \frac{\rho}{2}  \|S^\top (Ax-  b)\|^2_{G^{-1}} \nonumber \\
    & = \argmin_x \ \frac{1}{2} \|x-x_k\|^2_B + z_{k}^\top  S^\top (Ax-b) + \frac{\rho}{2}  \|S^\top (Ax- b)\|^2_{G^{-1}} \nonumber \\
    & \overset{\textbf{def}}{:=} \argmin_x \ \mathcal{L}(x,z_k, \rho).
\end{align} 
Later, we will derive the closed form of the update formulas provided in \eqref{als:ls:003} and \eqref{als:ls:004}. Before we delve into the details, first let us define the following matrices.
\begin{definition}
\label{als:ls:d:2}
We use throughout the following definitions
\begin{align*}
  & H := \left(\frac{1}{\rho} G + S^\top AB^{-1}A^\top S \right)^{-1},  \   Z := A^\top S H S^\top A, \\
  & F := \left(\frac{1}{\rho} B + A^\top SG^{-1}S^\top A \right)^{-1} , \ \   W := S^\top AF A^\top S.
\end{align*}
\end{definition}

\begin{remark}
The main difference between the PS and ALS method is that in the ALS method we maintain a dual update along with the primal update. Our framework can be sought as augmented Lagrangian method  \footnote{This is also known as the \textit{Method of Multipliers} \cite{Hestenes1969,powell,BERTSEKAS}.} applied to the optimization problem \eqref{als:ls:1}. 

\end{remark}


 \begin{algorithm}
\caption{$x_{K+1} = \textbf{ALS}(A,b, B, G, c, \mathcal{D}, K)$}
\label{als:ls:alg:als}
\begin{algorithmic}
\STATE{Choose initial point $x_0 \in \R^n, \ z_0 \in \R^p, \ \rho_0 \in \R$}
\WHILE{$k \leq K$}
\STATE{Draw a fresh sample $S \sim \mathcal{D}$ and update
\begin{align*}
    & z_{k+1} = \left(\frac{1}{\rho_k} G + S^\top AB^{-1}A^\top S \right)^{-1} [S^\top A x_k - S^\top b + \frac{1}{\rho_k} G  z_k]; \\
    & x_{k+1} = x_k -  B^{-1}  A^\top  S z_{k+1}; \\
    & \rho_{k+1} = c \rho_k;
\end{align*}
$k \leftarrow k+1$;}
\ENDWHILE
\end{algorithmic}
\end{algorithm}

\begin{lemma}
\label{als:ls:lem:1}
The updates given in \eqref{als:ls:003} and \eqref{als:ls:004} are equivalent to
\begin{align*}
 z_{k+1} = H \left(S^\top Ax_k-S^\top b +\frac{1}{\rho} Gz_k\right), \quad    x_{k+1} =  x_k - B^{-1}  A^\top  S z_{k+1}.
\end{align*}
\end{lemma}

\begin{proof}
From the definition, we have
\begin{align*}
    & x_{k+1} = \argmin_{x} \mathcal{L} (x,z_k, \rho) , \quad z_{k+1} = z_k + \rho G^{-1} S^\top  (Ax_{k+1}-b).
\end{align*}
Setting $\frac{\partial \mathcal{L} (x,z_k, \rho)}{\partial x} = 0$, we get
\begin{align*}
    0 = \frac{\partial \mathcal{L} (x,z_k, \rho)}{\partial x} = B(x-x_k)+ A^\top  S z_k + \rho  A^\top S G^{-1}(S^\top Ax-S^\top b).
\end{align*}
Rearranging the above along with the definition of $F$, we get the following:
\begin{align*}
    \rho F^{-1} x = \rho \left(\frac{1}{\rho} B + A^\top SG^{-1}S^\top A \right)^{-1} x = Bx_k - A^\top  S z_k + \rho  A^\top S G^{-1} S^\top b.
\end{align*}
Since, $\rho > 0$ we get 
\begin{align}
\label{als:ls:13}
    x_{k+1} = \frac{1}{\rho} F \left(Bx_k - A^\top  S z_k + \rho  A^\top S G^{-1}S^\top b\right).
\end{align}
Now, we will simplify each term of \eqref{als:ls:13} using Lemma \ref{als:ls:lem:13}. The first term can be simplified as
\begin{align}
\label{als:ls:14}
 \frac{1}{\rho} F Bx_k  & = \left(B^{-1}-B^{-1} Z B^{-1}\right) Bx_k  = x_k - B^{-1} Z x_k.
\end{align}
Moreover, using Lemma \ref{als:ls:lem:13} we can derive the following identity:
\begin{align}
\label{als:ls:15}
  Z  B^{-1} A^\top S  & = A^\top S H S^\top A B^{-1} A^\top S  = A^\top S H (H^{-1}-\frac{1}{\rho} G) \nonumber \\
  & = A^\top S -\frac{1}{\rho} A^\top S HG.  
\end{align}
Then, the second term of \eqref{als:ls:13} can be written as
\begin{align}
\label{als:ls:16}
  - \frac{1}{\rho} F A^\top  S z_k  & = - B^{-1} A^\top S z_k + B^{-1}  Z B^{-1}  A^\top S z_k = -  \frac{1}{\rho}  B^{-1} A^\top S H G z_k,
\end{align}
here, we used the identity from \eqref{als:ls:15}. Similarly, using \eqref{als:ls:15} the third term of \eqref{als:ls:13} can be simplified as follows:
\begin{align}
\label{als:ls:17}
   F  A^\top S G^{-1}S^\top b & = \rho (B^{-1}-B^{-1} Z B^{-1})  A^\top S G^{-1}S^\top b  = B^{-1}  A^\top S H S^\top b.
\end{align}
Substituting \eqref{als:ls:14}, \eqref{als:ls:16} and \eqref{als:ls:17} in \eqref{als:ls:13}, we get the expression of $x_{k+1}$. To prove the second expression for $x_{k+1}$, let's denote $u_k = S^\top Ax_k-S^\top b +\frac{1}{\rho} G z_k$. Then, we have,
\begin{align}
\label{als:ls:18}
    G^{-1} S^\top  (Ax_{k+1}-b) = G^{-1}S^\top  (Ax_{k}-b) -G^{-1} GS^\top A B^{-1} A^\top S H u_k.
\end{align}
Now, using the third identity of Lemma \ref{als:ls:lem:13}, we can verify the following identity:
\begin{align}
  \label{als:ls:19}
  S^\top A B^{-1} A^\top S H = I - \frac{1}{\rho} G H.
\end{align}
Using the expression of \eqref{als:ls:19} in \eqref{als:ls:18}, we get
\begin{align}
   G^{-1} S^\top  (Ax_{k+1}-b) = - \frac{1}{\rho} z_k + \frac{1}{\rho} H u_k. 
\end{align}
Now, using the definition of $z_{k+1}$, we get
\begin{align*}
    z_{k+1} = z_k + \rho G^{-1} S^\top  (Ax_{k+1}-b) = H u_k.
\end{align*}
Substituting the above into the expression of $x_{k+1}$, we get the required identity.
\end{proof}

\begin{remark}
\label{als:ls:rem:1}
In Algorithm \ref{als:ls:alg:als} we provide an adaptive version of the ALS method. In the adaptive version, the penalty parameter $\rho$ is updated gradually so that $\rho_{k}$ becomes very large after finitely many iterations. To achieve this, we set a multiplier $c>1$ and update the penalty parameter using 
\begin{align*}
\rho_{k+1} = c \rho_k.
\end{align*}
This scheme improves the computational efficiency of the proposed ALS method.
This type of iterative approach is frequently used in the classical augmented Lagrangian framework \cite{Hestenes1969,BERTSEKAS}. The proposed penalty and augmented type algorithms have been proposed recently in the context of Kaczmarz method for solving linear systems and linear feasibility problems \cite{als:kac} as well as Newton's method for solving unconstrained minimization \cite{als:newton}.
\end{remark}

In the next lemma, we will show that in the limiting case, the PS and ALS methods resolve into the SP method.

\begin{lemma}
\label{als:ls:lem:2}
The PS and ALS methods resolve into the SP method for $\rho \rightarrow \infty$, i.e., $\lim_{\rho \rightarrow \infty} \text{PS} \equiv \lim_{\rho \rightarrow \infty} \text{ALS} \equiv \text{SP}$.
\end{lemma}

\begin{proof}
Denote, $Y = B^{-\frac{1}{2}} A^\top S G^{-\frac{1}{2}}$ and take, $x^*$ such that $b =Ax^*$ holds. Then, from the update formula of the PS method, we have
\begin{align}
\label{als:ls:20}
 \lim_{\rho \rightarrow \infty} x_{k+1} & = x_k - B^{-1} \lim_{\rho \rightarrow \infty} A^\top S H S^\top A(x_k-x^*) \nonumber \\
 & = x_k - B^{-\frac{1}{2}} Y \lim_{\rho \rightarrow \infty} \left(\frac{1}{\rho}I+ Y^\top  Y\right)^{-1} Y^\top  B^{\frac{1}{2}} (x_k-x^*) \nonumber \\
 & = x_k - B^{-\frac{1}{2}} Y Y^{\dagger} B^{\frac{1}{2}} (x_k-x^*) \nonumber \\
 & = x_k - B^{-\frac{1}{2}} Y (Y^\top Y)^{\dagger} Y^\top  B^{\frac{1}{2}} (x_k-x^*) \nonumber \\
 & = x_k - B^{-1} A^\top S G^{-\frac{1}{2}}\left(G^{-\frac{1}{2}}S^\top AB^{-1}A^\top S G^{-\frac{1}{2}} \right)^{\dagger} G^{-\frac{1}{2}} S^\top A(x_k-x^*) \nonumber \\
 & = x_k - B^{-1} A^\top S (S^\top AB^{-1}A^\top S)^{\dagger}S^\top (Ax_k-b).
\end{align}
here, we used the following well-known identities of matrix pseudo-inverse (Theorem 4.3, \cite{Barata2012}):
\begin{align*}
 Y^{\dagger} =  \lim_{\rho \rightarrow \infty} (\frac{1}{\rho}I+ Y^\top  Y)^{-1} Y^\top   =  (Y^\top Y)^{\dagger} Y^\top .
\end{align*}
Denote, $u_k = S^\top Ax_k-S^\top b +\frac{1}{\rho} G z_k$. Since, $\lim_{\rho \rightarrow \infty} u_k = S^\top Ax_k-S^\top b$, using the derivation of \eqref{als:ls:20} for the ALS update we get
\begin{align}
\label{als:ls:21}
 \lim_{\rho \rightarrow \infty} x_{k+1} = x_k - B^{-1} A^\top S (S^\top AB^{-1}A^\top S)^{\dagger}S^\top (Ax_k-b).
\end{align}
The update formulas derived in \eqref{als:ls:20} and \eqref{als:ls:21} are the same and are precisely the update formula of the SP method (see \cite{gower:2015, richtrik2017stochastic}). Moreover, for the ALS method, in the limit $\rho \rightarrow \infty$, the expression for $z_{k+1}$ trivially holds, i.e.,
\begin{align*}
  0 = \lim_{\rho \rightarrow \infty} \left( \frac{1}{\rho}z_{k+1}-\frac{1}{\rho}z_k \right)= \lim_{\rho \rightarrow \infty}  G^{-1} S^\top  (Ax_{k+1}-b) = 0.
\end{align*}
In the last line we used the fact $S^\top (Ax_{k+1}-b) = 0$ (identity \eqref{als:ls:21} can be used to verify this). This proves the lemma.
\end{proof}

Before we delved into the technical details of the convergence theory, first let us specify the solution of problem \eqref{als:ls:1}. We derive closed-form solution of the problem \eqref{als:ls:1} as follows:
\begin{align}
\label{als:ls:22}
    x^*  = \argmin_{x \in \mathcal{X}} \|x_0-x\|_B^2  = x_0 - B^{-1} A^\top  \left(AB^{-1} A^\top \right)^{\dagger} (Ax_0-b),
\end{align}
here, $ \mathcal{X} = \{x \ | \ Ax = b\}$ and $x^*$ can be interpreted as the projection of the initial point $x_0$ onto the feasible region $ \mathcal{X}$ in the $B-$ norm. Note that, the choice of $x^*$ allows the vector $B^{\frac{1}{2}}(x_k-x^*)$ to remain in the range space of $B^{-1} A^\top $, i.e., $B^{\frac{1}{2}}(x_k-x^*) \in \textbf{Range}(B^{-1} A^\top )$. We will use this relation throughout the paper. Before that, let us define the following function:
\begin{definition}
\label{als:ls:d:3}
For any $x \in \R^n$ and $z \in \R^p, \ \rho > 0$, define the functions $f_S(x,z, \rho)$ and $f(x,z, \rho)$ as follows:
\begin{align*}
f_S(x, z, \rho) = \frac{1}{2} \|x-x^*\|^2_{Z} + \frac{1}{2\rho} \|z\|_{W}^2, \quad  f(x, z, \rho)  = \E_{\mathcal{D}}[f_S(x, z, \rho)].
\end{align*}
\end{definition}
Special cases of the above functions were first studied in \cite{richtrik2017stochastic} in the context of the SP method. The above functions can be sought as extensions of the the following functions defined in \cite{richtrik2017stochastic}:
\begin{align*}
f_S(x) = \frac{1}{2} \|x-x^*\|^2_{Z} , \quad  f(x)  = \E_{\mathcal{D}}[f_S(x)].
\end{align*}
Indeed, we can show that in the limiting case the functions defined in this work resolves into the functions defined in \cite{richtrik2017stochastic}, i.e.,
\begin{align*}
\lim_{\rho \rightarrow \infty} f_S(x, z, \rho) = f_S(x), \quad  \lim_{\rho \rightarrow \infty} f(x, z, \rho)  = f(x).
\end{align*}

Most of the convergence results derived in this work depend on the characteristics of $ f(x, z, \rho)$ (see Lemma \ref{als:ls:lem:6}). Until otherwise mentioned, throughout the paper, we assume sketch matrix $S$ satisfies the following assumption:

\begin{assumption}
\label{als:ls:a:1}
Let $S \sim \mathcal{D}$ be such that the following holds:
\begin{align*}
    \textbf{Null}\left(\E_{\mathcal{D}}[SHS^\top ]\right) \subset \textbf{Null}\left(A^\top  \right), \quad \E_{\mathcal{D}}[W] \succ \mathbf{0}.
\end{align*}
\end{assumption}
The first part of Assumption~\ref{als:ls:a:1} was introduced in \cite{richtrik2017stochastic} in the context of the SP method. Later, we discuss specific variations of Assumption \ref{als:ls:a:1} that will be used in this work. Assumption \ref{als:ls:a:1} allows us to derive the following equivalent condition between range spaces.

\begin{lemma}
\label{als:ls:lem:3}
We have the following:
\begin{align*}
  B^{\frac{1}{2}}(x_k-x^*) \in \textbf{Range}(B^{-1} A^\top ) =  \textbf{Range} \left(\E_{\mathcal{D}}[B^{-\frac{1}{2}} Z B^{-\frac{1}{2}}]\right).
\end{align*}
\end{lemma}

\begin{proof}
Form our choice of $x^*$, we have $B^{\frac{1}{2}}(x_k-x^*) \in \textbf{Range}(B^{-1} A^\top ) $. That means we just need to show that $\textbf{Range}(B^{-1} A^\top ) $ and $\textbf{Range} \left(\E_{\mathcal{D}}[B^{-\frac{1}{2}} Z B^{-\frac{1}{2}}]\right) $ are equivalent. Note that we have,
\begin{align*}
    \textbf{Null} & (\E_{\mathcal{D}}[B^{-\frac{1}{2}} Z B^{-\frac{1}{2}}]) = \textbf{Null} \left(B^{-\frac{1}{2}}A^\top  \E_{\mathcal{D}}[S H S^\top ] A B^{-\frac{1}{2}}\right) \\
    & =  \textbf{Null} \left(A^\top  \E_{\mathcal{D}}[S H S^\top ]  A B^{-\frac{1}{2}}\right) =  \textbf{Null} \left( A B^{-\frac{1}{2}}\right),
\end{align*}
where, the last equality follows from Lemma \ref{als:ls:lem:12} with the choice $P = \E_{\mathcal{D}}[S H S^\top ]$ and $ M = A$. Taking orthogonal complement, we get the required relation.
\end{proof}

Lemma \ref{als:ls:lem:3} is the so-called exactness condition. It was first introduced in \cite{richtrik2017stochastic} in the context of the SP method. This lemma will later be used to derive meaningful bounds that are crucial for the convergence analysis of the proposed algorithms (see Lemma \ref{als:ls:lem:6}). Using the above construction, we propose the following equivalent reformulation of the original problem.

\paragraph{Penalty \& Augmented Stochastic Reformulation.} Consider the following stochastic optimization problems:
\begin{align}
 & \textbf{PSR:} \quad (x^*) = \argmin_{x}  f(x, 0, \rho) = \frac{1}{2}\|x-x^*\|^2_{\E_{\mathcal{D}}[Z]},  \label{als:ls:23} \\
  & \textbf{ASR:} \quad (x^*,z^*) = \argmin_{x, z}  f(x, z, \rho) = \frac{1}{2}\|x-x^*\|^2_{\E_{\mathcal{D}}[Z]}+ \frac{1}{2\rho} \|z\|^2_{\E_{\mathcal{D}}[W]}.  \label{als:ls:24}
\end{align}
Note that, if we assume $\E_{\mathcal{D}}[Z] \succ \mathbf{0} $ and $\E_{\mathcal{D}}[W] \succ \mathbf{0} $, then we can easily show that the above optimization problems are equivalent to solving the linear system. However, in Assumption \ref{als:ls:a:1} we impose a weaker condition on the sketching matrix $S$. Under these circumstances, we will show that indeed the reformulations  \eqref{als:ls:23} and \eqref{als:ls:24} are equivalent to \eqref{als:ls:1} (see the following lemma).

\begin{lemma}
\label{als:ls:lem:4}
If Assumption \ref{als:ls:a:1} holds and there exists $\mu_{\mathcal{D}} > 0 $ such that
\begin{align*}
f(x,z, \rho) =  \|x-x^*\|^2_{\E_{\mathcal{D}}[Z]}+ \frac{1}{\rho} \|z\|^2_{\E_{\mathcal{D}}[W]}  \geq \mu_{\mathcal{D}} \left(\|x-x^*\|^2_B + \frac{1}{\rho} \|z\|^2_{G}\right)
\end{align*}
holds for all $x \in \R^n, z \in \R^p$, then the penalty and augmented reformulations provided in \eqref{als:ls:23} and \eqref{als:ls:24} are equivalent to the linear system \eqref{als:ls:1}.
\end{lemma}

\begin{proof}
We will prove the equivalence for the second reformulation provided in \eqref{als:ls:24}. The result then follows since \eqref{als:ls:23} is a special case of \eqref{als:ls:24}. First, let us rewrite the original linear system as follows:
\begin{align*}
    Ax - b = z, \quad z = 0
\end{align*}
We can easily verify that the above system is indeed equivalent to the linear system $Ax = b$. Assume, $\bar{x}, \bar{z}$ are the solutions of the above system. Now, substituting the values $\bar{x} = x^*, \bar{z} = 0$ in \eqref{als:ls:24}, we have
\begin{align*}
    f(\bar{x}, \bar{z}, \rho) = \|\bar{x}-x^*\|_{\E_{\mathcal{D}}[Z]} + \frac{1}{\rho} \|\bar{z}\|^2_{\E_{\mathcal{D}}[W]} = 0
\end{align*}
Since, the optimal objective value of optimization problem is $0$, we can say that the pair $(\bar{x}, \bar{z})$ solves the optimization problem \eqref{als:ls:24}. 

On the other hand, assume the pair $(\Tilde{x},\Tilde{z})$ is the optimal solution of problem \eqref{als:ls:24}. This implies the following:
\begin{align*}
& 0 = f(\Tilde{x},\Tilde{z}, \rho) =  \|\Tilde{x}-x^*\|^2_{\E_{\mathcal{D}}[Z]}+ \frac{1}{\rho} \|\Tilde{z}\|^2_{\E_{\mathcal{D}}[W]}  \geq \mu_{\mathcal{D}} [\|\Tilde{x}-x^*\|^2_B + \frac{1}{\rho} \|\Tilde{z}\|^2_{G}]. \\
\Longleftrightarrow \quad & 0  \geq \mu_{\mathcal{D}} [\|\Tilde{x}-x^*\|^2_B + \frac{1}{\rho} \|\Tilde{z}\|^2_{G}] \\
\Longleftrightarrow \quad &  \Tilde{x}-x^* = 0, \quad \Tilde{z} = 0
\end{align*}
This implies that the pair $(\Tilde{x},\Tilde{z})$ solves the linear system. This proves that the linear system and the optimization problem \eqref{als:ls:24} are equivalent.
\end{proof}

Later, we will show that for special sketching techniques constant $\mu_{\mathcal{D}} > 0$ exists. Note that, Lemma \ref{als:ls:lem:4} allows us to solve the penalty/augmented sketching problems instead of the original problem. In a similar fashion, if we take $\rho \rightarrow \infty$ we get the equivalent optimization problem defined in \cite{richtrik2017stochastic}. However, the proposed ALS method chooses a reasonable $\rho$ that generates the acceleration in the algorithmic performance. Similarly, for the adaptive ALS method, the sequence $\rho_k$ becomes sufficiently large after finitely many iterations that provide a good approximation for infinite $\rho$ value. Furthermore, later we will discuss the above exactness condition in the context of several sampling strategies.

\section{Convergence Theory}
\label{als:ls:sec:conv}
We now present the convergence results for the proposed PS and ALS methods. Our goal is to show that the primal sequence $x_{k}$ generated by both the PS and ALS method converges. Furthermore, the dual sequence $z_k$ generated by the ALS algorithm converges. In brief, we will show the following:
\begin{align*}
& \E[\|x_{k}-x^*\|^2_B] \rightarrow 0, \quad   \E[\|z_k\|^2_{G}] \rightarrow 0, \\
& \E[\|x_{k}-x^*\|^2_{\E_{\mathcal{D}}[Z]}] \rightarrow 0, \quad \E[\|z_k\|^2_{\E_{\mathcal{D}}[W]}] \rightarrow 0.
\end{align*}
We start by showing the following result regarding the eigenvalues of matrix $B^{-\frac{1}{2}} Z B^{-\frac{1}{2}}$ and $G^{-\frac{1}{2}} W G^{-\frac{1}{2}}$. 

\begin{lemma}
\label{als:ls:lem:5}
For any $S \sim \mathcal{D}$, we have the following:
\begin{align*}
    \lambda_{i} \left(B^{-\frac{1}{2}} Z B^{-\frac{1}{2}} \right)  = \lambda_{i} \left(B^{-\frac{1}{2}} Z B^{-\frac{1}{2}}  \right) =  \lambda_{i} \left(G^{\frac{1}{2}} W G^{\frac{1}{2}}  \right) =  \lambda_{i} \left(G^{-\frac{1}{2}} W G^{-\frac{1}{2}} \right),
\end{align*}
here, $\lambda_i(\cdot)$ denotes the $i^{th}$ eigenvalue.
\end{lemma}

\begin{proof}
Take, $Y = B^{-\frac{1}{2}} A^\top S G^{-\frac{1}{2}}$. Assume, matrix $Y$ has the following singular value decomposition:
\begin{align*}
    Y = U \Sigma V^\top , \quad Y^\top   = V \Sigma^\top  U^\top ,
\end{align*}
where, $U$ and $V$ are the corresponding orthogonal matrices. Using the above decomposition, we have
\begin{align*}
    &  B^{-\frac{1}{2}} Z B^{-\frac{1}{2}} = Y \left(\frac{1}{\rho} I + Y^\top Y \right)^{-1} Y^\top   = U \Sigma  \left(\frac{1}{\rho} I + \Sigma^\top \Sigma \right)^{-1} \Sigma^\top  U^\top ,  \\
    &  G^{-\frac{1}{2}} W G^{-\frac{1}{2}} = Y^\top  \left(\frac{1}{\rho} I + YY^\top  \right)^{-1} Y = V \Sigma^\top   \left(\frac{1}{\rho} I + \Sigma\Sigma^\top  \right)^{-1} \Sigma V^\top .
\end{align*}
Since, both $\Sigma  \left(\frac{1}{\rho} I + \Sigma^\top \Sigma \right)^{-1} \Sigma^\top $ and $\Sigma^\top   \left(\frac{1}{\rho} I + \Sigma\Sigma^\top  \right)^{-1} \Sigma$ matrices are diagonal matrices with identical nonzero entries, we get the required result.
\end{proof}
The above lemma allows us to obtain a reasonable convergence rate bound. In general, the convergence rate of the proposed methods depends on the minimum positive eigenvalues of $\E_{\mathcal{D}}[B^{-\frac{1}{2}} Z B^{-\frac{1}{2}}]$ and $\E_{\mathcal{D}}[G^{-\frac{1}{2}} W G^{-\frac{1}{2}}]$. In the following, we define function $\mathcal{V}(x,z, \rho) $.

\begin{definition}
\label{als:ls:d:4}
For $k \geq 1$, assume $x_k$ and $z_k$ are random iterates of the ALS method. Then, define the following function:
\begin{align}
\label{als:ls:25}
\mathcal{V}(x_k,z_k, \rho)  = \|x_k-x^*\|^2_B + \frac{1}{\rho} \|z_k\|_{G}^2.
\end{align}
\end{definition}
Furthermore, we define $\mathcal{V}(\bar{x}_k,\bar{z}_k, \rho) $ as the corresponding function with respect to the Cesaro averages, $\bar{x}_k$ and $\bar{z}_k$, i.e., $\bar{x}_k = \sum \nolimits_{l=0}^{k-1} x_l$ and $\bar{z}_k = \sum \nolimits_{l=0}^{k-1} z_l$. Throughout the convergence analysis, our goal is to bound the growth of functions $\E[\mathcal{V}(x_k,z_k, \rho) ]$ and $\E[f(x_k,z_k, \rho) ]$. indeed, we will show that for the ALS method, $\E[\mathcal{V}(x_k,z_k, \rho) ]$ is a Lyapunov function.

\begin{lemma}
\label{als:ls:lem:6}
\sloppy Assume Assumption \ref{als:ls:a:1} holds. Denote, $\mu = \min\{\lambda_{\min}^+(\E_{\mathcal{D}}[B^{-\frac{1}{2}} Z B^{-\frac{1}{2}}])$, $\lambda_{\min}^+(\E_{\mathcal{D}}[G^{-\frac{1}{2}} W G^{-\frac{1}{2}}])\}$ and $\eta =  \max \{\lambda_{\max}(\E_{\mathcal{D}}[B^{-\frac{1}{2}} Z B^{-\frac{1}{2}}])$, $\lambda_{\max}(\E_{\mathcal{D}}[G^{-\frac{1}{2}} W G^{-\frac{1}{2}}])\}$. Then, if $x_k, \ z_k$ are the random iterates of the ALS algorithm, we have the following:
\begin{align*}
&  \E_{\mathcal{D}}[\lambda_{\min}^+(G^{-\frac{1}{2}} W G^{-\frac{1}{2}})] \ \mathcal{V}(x_k,z_k, \rho)  \ \leq \ \mu \mathcal{V}(x_k,z_k, \rho)   \ \leq \  2 f(x_k, z_k, \rho),  \\
& 2f(x_k, z_k, \rho) \ \leq \ \eta \mathcal{V}(x_k,z_k, \rho)  \ \leq \ \E_{\mathcal{D}}[\lambda_{\max}(G^{-\frac{1}{2}} W G^{-\frac{1}{2}})] \ \mathcal{V}(x_k,z_k, \rho).
\end{align*}
\end{lemma}

\begin{proof}
From, the definition, we have
\begin{align}
    \label{als:ls:26}
    2f & (x_k, z_k, \rho)    = \|x_k-x^*\|^2_{\E_{\mathcal{D}}[Z]} + \frac{1}{\rho} \|z_k\|_{\E_{\mathcal{D}}[W]}^2  \nonumber \\
    & \leq \lambda_{\max}(\E_{\mathcal{D}}[B^{-\frac{1}{2}} Z B^{-\frac{1}{2}}] )\|x_k-x^*\|^2_B + \frac{1}{\rho} \lambda_{\max}(\E_{\mathcal{D}}[G^{-\frac{1}{2}} W G^{-\frac{1}{2}}])\} \|z_k\|^2_{G} \nonumber \\
    & \leq \eta \ \mathcal{V}(x_k,z_k, \rho),
\end{align}
here, we used the relations $\E_{\mathcal{D}}[Z] \preceq \lambda_{\max}(\E_{\mathcal{D}}[B^{-\frac{1}{2}} Z B^{-\frac{1}{2}}]) B$ and \\ $\E_{\mathcal{D}}[W] \preceq \lambda_{\max}(\E_{\mathcal{D}}[G^{-\frac{1}{2}} W G^{-\frac{1}{2}}]) G$. Since the function $\lambda_{\max}$ is convex over symmetric matrices the upper bounds follows from \eqref{als:ls:26} along with Jensens' inequality. From Lemma \ref{als:ls:lem:3} we have $B^{\frac{1}{2}}(x_k-x^*) \in \textbf{Range}(\E_{\mathcal{D}}[B^{-\frac{1}{2}} Z B^{-\frac{1}{2}}])$ and $\E_{\mathcal{D}}[G^{-\frac{1}{2}} W G^{-\frac{1}{2}}] \succ 0$ holds from our assumption. Now, using the Courant-Fischer theorem (Theorem 8.1.2, \cite{golub2013}), we have
\begin{align}
    \label{als:ls:27}
     2f & (x_k, z_k, \rho)   = \|x_k-x^*\|^2_{\E_{\mathcal{D}}[Z]} + \frac{1}{\rho} \|z_k\|_{\E_{\mathcal{D}}[W]}^2  \nonumber \\
    & \geq \lambda_{\min}^+(\E_{\mathcal{D}}[B^{-\frac{1}{2}} Z B^{-\frac{1}{2}}]) \|x_k-x^*\|^2_B + \frac{1}{\rho} \lambda_{\min}^+(\E_{\mathcal{D}}[G^{-\frac{1}{2}} W G^{-\frac{1}{2}}]) \|z_k\|^2_{G} \nonumber \\
    & \geq \mu \ \mathcal{V}(x_k,z_k, \rho).
\end{align}
Combining \eqref{als:ls:26} and \eqref{als:ls:27}, we get the required result.
\end{proof}

The above lemma bounds the sequence $f(x_{k},z_{k},\rho)$ in terms of the Lyapunov function $\mathcal{V}(x_k,z_k, \rho) $. One can interpret these bounds as the respective strong convexity condition and the smoothness condition of the function $f(x_{k},z_{k}, \rho)$. Moreover, in the limit $\rho \rightarrow \infty$ the function $f(x_{k},z_{k}, \rho)$ resolves into the function $f(x_k)$ defined in \cite{richtrik2017stochastic}.

\begin{lemma}
\label{als:ls:lem:7}
Denote, $ Y = B^{-\frac{1}{2}} A^\top S G^{-\frac{1}{2}}$ and $\sigma = \max_{S \sim \mathcal{D}} \lambda_{\max}(YY^\top )$. Then, the following holds:
\begin{align*}
&  0 < \mu \leq \eta \leq \frac{\rho \sigma}{1+ \rho \sigma} < 1.
\end{align*}
\end{lemma}

\begin{proof}
For simplification, denote, $\bar{W} = G^{-\frac{1}{2}} W G^{-\frac{1}{2}}$. Using the definition, we have
\begin{align}
\label{als:ls:28}
    \bar{W}^2 & = G^{-\frac{1}{2}} W G^{-1} W G^{-\frac{1}{2}}  = G^{-\frac{1}{2}} S^\top  A F A^\top  AG^{-1} S^\top A F A^\top S G^{-\frac{1}{2}} \nonumber \\
    & = \bar{W} - \frac{1}{\rho} G^{-\frac{1}{2}} S^\top  A F B F A^\top S G^{-\frac{1}{2}} \leq \bar{W} - \frac{\lambda_{\min}\left(F^{\frac{1}{2}}B F^{\frac{1}{2}} \right)}{\rho}  \bar{W}.
\end{align}
Now, we can simplify the following:
\begin{align}
\label{als:ls:29}
  \lambda_{\min}\left(F^{\frac{1}{2}}B F^{\frac{1}{2}} \right) & = \lambda_{\min}\left( \left(\frac{1}{\rho} + YY^\top \right)^{-1} \right)    = \frac{1}{\lambda_{\max}\left( \frac{1}{\rho} + YY^\top \right)} \nonumber \\
  & \geq \frac{\rho}{1+ \rho \lambda_{\max}(YY^\top )} \geq \frac{\rho}{1+ \rho \sigma}.
\end{align}
Substituting the above in \eqref{als:ls:28}, we have
\begin{align}
\label{als:ls:30}
    \bar{W}^2 & \leq \bar{W} - \frac{\lambda_{\min}\left(F^{\frac{1}{2}}B F^{\frac{1}{2}} \right)}{\rho}  \bar{W} \leq \frac{\rho \sigma}{1+ \rho \sigma}  \bar{W}.
\end{align}
Now, for any vector $y \neq 0 $, we have
\begin{align*}
  \|\E_{\mathcal{D}}[\bar{W}]y\|^2  =  y^\top  \left(\E_{\mathcal{D}}[\bar{W}]\right)^2 y & \leq  y^\top  \E_{\mathcal{D}}[\bar{W}^2] y \\
  & \leq \frac{\rho \sigma \ y^\top  \E_{\mathcal{D}}[\bar{W}] y}{1+ \rho \sigma}   \leq \frac{\rho \sigma \ \|\E_{\mathcal{D}}[\bar{W}]y\| \|y\|}{1+ \rho \sigma}.
\end{align*}
Simplifying further, we have 
\begin{align*}
    \lambda_{\max}({\E_{\mathcal{D}}[\bar{W}]}) = \max_{y \neq 0} \frac{\|\E_{\mathcal{D}}[\bar{W}]y\|}{\|y\|} \leq \frac{\rho \sigma}{1+ \rho \sigma} < 1.
\end{align*}
Since, $ \mu \leq \lambda_{\max}({\E_{\mathcal{D}}[\bar{W}]})$, we have the required result.
\end{proof}

In the following theorems, we provide global linear convergence of the sequences $x_k$ and $z_k$ generated by both the PS and ALS methods.

\begin{theorem}
\label{als:ls:thm:1}
If $x_k$ is the random iterate of the PS algorithm with fixed penalty $\rho$, then the sequence $x_k$ converges and the following results hold:
\begin{align*}
   & \E[\|x_{k+1}-x^*\|^2_B]  \leq \left(1-\mu_{\min}-\frac{\mu_0}{\rho}\right)^{k+1} \ \E[\|x_{0}-x^*\|^2_B],  \\
   & \E[\|x_{k+1}-x^*\|^2_{\E[Z]} ] \leq \eta_{\max} \left(1-\mu_{\min}-\frac{\mu_0}{\rho}\right)^{k+1} \ \E[\|x_{0}-x^*\|^2_B].
\end{align*}
Where, $\mu_{\min} = \lambda_{\min}^+(\E_{\mathcal{D}}[B^{-\frac{1}{2}} Z B^{-\frac{1}{2}}])$, \  $\eta_{\max} = \lambda_{\max}(\E_{\mathcal{D}}[B^{-\frac{1}{2}} Z B^{-\frac{1}{2}}])$, and \\ $\mu_0 = \lambda_{\min}( B^{-\frac{1}{2}}\E_{\mathcal{D}}[A^\top S HGHS^\top A] B^{-\frac{1}{2}} )$.
\end{theorem}

\begin{proof}
From the update formula of the PS method we have,
\begin{align}
\label{als:ls:31}
    \|x_{k+1}-x^*\|^2_B  & = \|(I- B^{-1}Z)(x_k-x^*)\|^2_B \nonumber \\
    & =(x_k-x^*)^\top  [B-2Z+ZB^{-1}Z] (x_k-x^*) \nonumber \\
  & = (x_k-x^*)^\top  [B-Z-\frac{1}{\rho} A^\top S HGHS^\top A] (x_k-x^*) \nonumber \\
  & = \|x_k-x^*\|^2_B - \|x_k-x^*\|^2_Z - \frac{1}{\rho} \|HS^\top A (x_k-x^*)\|^2_{G},
\end{align}
here, we used the identity $ZB^{-1}Z = Z-\frac{1}{\rho} A^\top S HG HS^\top A$ from Lemma \ref{als:ls:lem:13}. Taking expectation in \eqref{als:ls:31} and simplifying, we have
\begin{align}
\label{als:ls:32}
     \E_{\mathcal{D}}[\|x_{k+1} & -x^*\|^2_B ]  \nonumber \\
     & =  \|x_k-x^*\|^2_B - \|x_{k}-x^*\|^2_{\E_{\mathcal{D}}[Z]} -   \frac{1}{\rho} \|x_{k}-x^*\|^2_{\E_{\mathcal{D}}[A^\top S HGHS^\top A]} \nonumber \\
    & \leq \  \|x_k-x^*\|^2_B - \mu_{\min} \|x_{k}-x^*\|^2_{B}- \frac{\mu_0}{\rho} \|x_{k}-x^*\|^2_{B} \nonumber  \\
    & \leq \ \left(1-\mu_{\min}-\frac{\mu_0}{\rho}\right) \|x_{k}-x^*\|^2_{B},
\end{align}
where we used that $$\E_{\mathcal{D}}[Z] \succeq \mu_{\min} B $$ and
$$\E_{\mathcal{D}}[A^\top S HGHS^\top A] \succeq \mu_{0}  $$ in the first and second inequality, respectively. 
Taking expectation again in \eqref{als:ls:32} along with the tower property of expectation and unrolling the recurrence, we get the first result of Theorem \ref{als:ls:thm:1}. For proving the second result, note that the following holds:
\begin{align*}
 \E[\|x_{k+1}-x^*\|^2_{\E_{\mathcal{D}}[Z]} ]  & \leq \ \lambda_{\max}(\E_{\mathcal{D}}[B^{-\frac{1}{2}} Z B^{-\frac{1}{2}}]) \ \E[ \|x_{k+1}-x^*\|^2_{B}]  \\
  & \leq \ \eta_{\max}  \left(1-\mu_{\min}-\frac{\mu_0}{\rho}\right)^{k+1}\ \E[\|x_{0}-x^*\|^2_B].
\end{align*} 
This proves the second part of Theorem \ref{als:ls:thm:1}.
\end{proof}

\begin{theorem}
\label{als:ls:thm:2}
Assume, $x_k, \ z_k$ are random iterates of the ALS algorithm with fixed penalty $\rho$, then the sequences $x_k, \ z_k$ converge and the following results hold:
\begin{align*}
   & \E[\mathcal{V}(x_{k+1},z_{k+1}, \rho) ] \leq (1-\mu)^{k+1} \ \E[\mathcal{V}(x_{0},z_{0}, \rho) ], \\
   & \E[f(x_{k+1},z_{k+1}, \rho)] \leq \frac{\eta}{2} (1-\mu)^{k+1} \E[\mathcal{V}(x_{0},z_{0}, \rho) ].
\end{align*}
\end{theorem}

\begin{proof}
For simplification, denote, $v_k = \frac{1}{\rho} G z_k$. Then, from the update formula of the ALS method we have,
\begin{align}
\label{als:ls:33}
    \|x_{k+1}-x^*\|^2_B  & = \|(I- B^{-1}Z)(x_k-x^*)- B^{-1}A^\top S H v_k\|^2_B \nonumber \\
    & = \|(I- B^{-1}Z)(x_k-x^*)\|^2_B + v_k^\top  HS^\top A B^{-1}A^\top S H v_k \nonumber \\
    & \quad \quad \quad \quad - 2 (x_k-x^*)^\top  [A^\top S H-ZB^{-1} A^\top S H] v_k.
\end{align}
The first term of \eqref{als:ls:33} can be simplified as follows:
\begin{align}
    \label{als:ls:34}
   \|(I- B^{-1}Z) & (x_k-  x^*)  \|^2_B  =(x_k-x^*)^\top  [B-2Z+ZB^{-1}Z] (x_k-x^*) \nonumber \\
  & = (x_k-x^*)^\top  [B-Z-\frac{1}{\rho} A^\top S HGHS^\top A] (x_k-x^*) \nonumber \\
  & = \|x_k-x^*\|^2_B - \|x_k-x^*\|^2_Z - \frac{1}{\rho} \|HS^\top A (x_k-x^*)\|^2_{G},
\end{align}
here, we used the identity $ZB^{-1}Z = Z-\frac{1}{\rho} A^\top S HG HS^\top A$ from Lemma \ref{als:ls:lem:13}. We can simplify the second term of \eqref{als:ls:33} as follows:
\begin{align}
\label{als:ls:35}
  v_k^\top  HS^\top A B^{-1}A^\top S H v_k  & =   v_k^\top  [H-\frac{1}{\rho} H G H] v_k  = \|v_k\|^2_H - \frac{1}{\rho} \|Hv_k\|_{G}^2.
\end{align}
Similarly, we can simplify the third term of \eqref{als:ls:33} as follows:
\begin{align}
    \label{als:ls:36}
    - 2 (x_k-x^*)^\top  & [A^\top S H-ZB^{-1} A^\top S H] v_k = -\frac{2}{\rho} (x_k-x^*)^\top   A^\top S HG H v_k.
\end{align}
Here, we used the identities from Lemma \ref{als:ls:lem:13}. Substituting the identities from \eqref{als:ls:34}-\eqref{als:ls:36} in \eqref{als:ls:33}, we get the following:
\begin{align}
\label{als:ls:37}
    & \|x_{k+1}-x^*\|^2_B \nonumber \\
    & = \|(I- B^{-1}Z)(x_k-x^*)\|^2_B + v_k^\top  HS^\top A B^{-1}A^\top S H v_k  \nonumber \\
    & \quad \quad \quad \quad - 2 (x_k-x^*)^\top  [A^\top S H-ZB^{-1} A^\top S H] v_k \nonumber \\
    & = \|x_k-x^*\|^2_B - \|x_k-x^*\|^2_Z - \frac{1}{\rho} \|HS^\top A (x_k-x^*)\|^2_{G} + \|v_k\|^2_H - \frac{1}{\rho} \|Hv_k\|_{G}^2 \nonumber \\
    & \quad \quad \quad  -\frac{2}{\rho} (x_k-x^*)^\top   A^\top S HG H v_k.
\end{align}
Now, using the dual update we have,
\begin{align}
   \label{als:ls:38} 
   \rho & \|v_{k+1}\|_{G^{-1}}^2 = \frac{1}{\rho} \|HS^\top A (x_k-x^*) + Hv_k\|_{G}^2 \nonumber \\
   & \quad \quad  \quad \quad  = \frac{1}{\rho} \|HS^\top A (x_k-x^*)\|_{G}^2 + \frac{1}{\rho} \|Hv_k\|_{G}^2 + \frac{2}{\rho}  (x_k-x^*)^\top A^\top S H G Hv_k.
\end{align}
Collecting terms from \eqref{als:ls:37} and \eqref{als:ls:38} and using the definition, we get
\begin{align}
\label{als:ls:39}
 \mathcal{V} & (x_{k+1},z_{k+1}, \rho)  = \|x_{k+1}-x^*\|^2_B + \rho\|v_{k+1}\|_{G^{-1}}^2 \nonumber \\
 & = \|x_k-x^*\|^2_B - \|x_k-x^*\|^2_Z - \cancel{\frac{1}{\rho} \|HS^\top A (x_k-x^*)\|^2_{G}} + \|v_k\|^2_H  \nonumber \\
   &  - \cancel{\frac{1}{\rho} \|Hv_k\|_{G}^2}   - \cancel{\frac{2}{\rho} (x_k-x^*)^\top   A^\top S HG H v_k} \nonumber \\
   & + \cancel{\frac{1}{\rho} \|HS^\top A (x_k-x^*)\|_{G}^2 }+ \cancel{\frac{1}{\rho} \|Hv_k\|_{G}^2} + \cancel{\frac{2}{\rho}  (x_k-x^*)^\top A^\top S H G Hv_k}  \nonumber \\
 & = \|x_{k}-x^*\|^2_B- \|x_{k}-x^*\|^2_Z + \|v_k\|^2_H.
\end{align}
Using the second identity of Lemma \ref{als:ls:lem:13}, we can simplify the following:
\begin{align}
\label{als:ls:40}
\|v_k\|^2_H = \frac{1}{\rho^2}\|z_k\|^2_{GHG}  = \frac{1}{\rho} \|z_k\|_{G} -   \frac{1}{\rho} \|z_k\|^2_{W}.
\end{align}
Using \eqref{als:ls:40} in \eqref{als:ls:39}, we have
\begin{align}
    \label{als:ls:41}
    \mathcal{V}(x_{k+1},z_{k+1}, \rho)  & = \mathcal{V}(x_{k},z_{k}, \rho) - \|x_{k}-x^*\|^2_Z -   \frac{1}{\rho} \|z_k\|^2_{W} \nonumber \\
    & = \mathcal{V}(x_{k},z_{k}, \rho)  - 2f_S(x_{k},z_{k}, \rho).
\end{align}
Taking expectation and using Lemma \ref{als:ls:lem:6}, we have
\begin{align}
     \E_{\mathcal{D}}[\mathcal{V}(x_{k+1},z_{k+1}, \rho) ]   & =  \ \mathcal{V}(x_{k},z_{k}, \rho) - 2f(x_{k},z_{k}, \rho)  \label{als:ls:42} \\
    & \leq \ (1-\mu) \ \mathcal{V}(x_k,z_k, \rho).      \label{als:ls:43}
\end{align}
Taking expectation again in \eqref{als:ls:43} along with the tower property of expectation and unrolling the recurrence, we get the first result of Theorem \ref{als:ls:thm:2}. For proving the second result, note that the following holds:
\begin{align*}
  2 \E[f(x_{k+1},z_{k+1} &, \rho)] = \E \left[\|x_{k+1}-x^*\|^2_{\E_{\mathcal{D}}[Z]} + \frac{1}{\rho} \|z_{k+1}\|_{\E_{\mathcal{D}}[W]}^2 \right] \\
  & \leq \ \lambda_{\max}(\E_{\mathcal{D}}[B^{-\frac{1}{2}} Z B^{-\frac{1}{2}}]) \ \E[\|x_{k+1}-x^*\|^2_{B}] \\
  & + \frac{\lambda_{\max}(\E_{\mathcal{D}}[G^{-\frac{1}{2}} W G^{-\frac{1}{2}}])}{\rho} \E[\|z_{k+1}\|_{G}^2] \\
  & \leq \ \eta \ \E[\mathcal{V}(x_{k+1},z_{k+1}, \rho)]  \leq \eta \ (1-\mu)^{k+1} \E[\mathcal{V}(x_{0},z_{0}, \rho) ].
\end{align*} 
This proves the second part of Theorem \ref{als:ls:thm:2}.
\end{proof}

By virtue of lemma \ref{als:ls:lem:7}, we can guarantee that the convergence rate obtained in the theorem \ref{als:ls:thm:2} is always bounded between one and zero, i.e., $0 < \frac{1}{1+\rho \sigma } \leq 1-\mu < 1$.

\begin{remark}
\label{als:ls:rem:2}
For simplicity, we proved Theorem \ref{als:ls:thm:2} with fixed $\rho$. The proof follows the same argument for the case of adaptive penalty parameter, i.e., $\rho_k$. Indeed, considering \eqref{als:ls:41}, we have the following:
\begin{align}
\label{als:ls:44}
    \mathcal{V}(x_{k+1},z_{k+1}, \rho_{k+1})  & = \mathcal{V}(x_{k},z_{k}, \rho_k) - \|x_{k}-x^*\|^2_Z  \nonumber \\
    & \quad \quad \quad \quad -   \frac{1}{\rho_k} \|z_k\|^2_{W} -  \frac{c-1}{c \rho_k} \|z_{k+1}\|_{G}^2.
\end{align}
With $c > 1$, we will always obtain a better convergence than the case with $c =1$. In that regard, one needs to obtain a reasonable lower bound of the quantity $\E[\|z_{k+1}\|^2_{G}]$ in terms of $\mathcal{V}(x_k,z_k)$. Moreover, the choice $c= 1$ resolves into the case with fixed penalty $\rho$.
\end{remark}

In the next theorem, we provide the convergence results of the average iterates $\bar{x}_k = \sum \nolimits_{l=0}^{k-1} x_l$ and $\bar{z}_k = \sum \nolimits_{l=0}^{k-1} z_l$ \footnote{This is widely known as Cesaro average.} generated by the ALS algorithm. We prove $\mathcal{O}(\frac{1}{k})$ convergence rate for the proposed ALS method.

\begin{theorem}
\label{als:ls:thm:3}
The average iterates $\bar{x}_k = \sum \nolimits_{l=0}^{k-1} x_l$ and $\bar{z}_k = \sum \nolimits_{l=0}^{k-1} z_l$ generated by the ALS method satisfy the following:
\begin{align*}
  \E[\mathcal{V}(\bar{x}_{k},\bar{z}_{k}, \rho)  ] = & \E[\|\bar{x}_k-x^*\|^2_{B}] + \frac{1}{\rho} \E[\|\bar{z}_k\|_{G}^2] \leq \frac{1}{\mu k} \mathcal{V}(x_{0},z_{0}, \rho),  \\
 \E[f(\bar{x}_{k},\bar{z}_{k}, \rho) ] =  & \E[\|\bar{x}_k-x^*\|^2_{\E_{\mathcal{D}}[Z]}] + \frac{1}{\rho} \E[\|\bar{z}_k\|_{\E_{\mathcal{D}}[W]}^2] \leq \frac{1}{k} \mathcal{V}(x_{0},z_{0}, \rho).
\end{align*}
\end{theorem}

\begin{proof}
Using \eqref{als:ls:42}, we have the following:
\begin{align}
    \label{als:ls:45}
    \sum \limits_{l =0}^{k-1} 2f(x_{l},z_{l} & , \rho)  = \sum \limits_{l =0}^{k-1} \|x_{l}-x^*\|^2_{\E_{\mathcal{D}}[Z]} + \frac{1}{\rho} \sum \limits_{l =0}^{k-1} \|z_l\|^2_{\E_{\mathcal{D}}[W]}  \nonumber \\
    & = \sum \limits_{l =0}^{k-1}  \E [\mathcal{V}(x_{l},z_{l}, \rho) ]- \E[\mathcal{V}(x_{l+1},z_{l+1}, \rho) ] \ \leq \ \mathcal{V}(x_{0},z_{0}, \rho).
\end{align}
Now, we deduce
\begin{align*}
   \E[\|\bar{x}_k-x^*\|^2_{B}] +  \frac{1}{\rho}  \E[\|\bar{z}_k\|_{G}^2]  & \leq  \E \left[\Big \| \frac{1}{k} \sum \limits_{l=0}^{k-1} \left(x_l-x^*\right)\Big \|^2_B\right] +  \frac{1}{\rho}\E \left[\Big \| \frac{1}{k} \sum \limits_{l=0}^{k-1} z_l\Big \|^2_{G}\right] \\
   & \leq  \E \left[\frac{1}{k} \sum \limits_{l=0}^{k-1}  \| x_l-x^* \|^2_B\right] + \frac{1}{\rho} \E \left[\frac{1}{k} \sum \limits_{l=0}^{k-1}  \| z_l \|^2_{G}\right] \\
   & = \frac{1}{k} \sum \limits_{l=0}^{k-1} \E[\mathcal{V}(x_{l},z_{l}, \rho) ] \leq \frac{1}{\mu k} \sum \limits_{l=0}^{k-1} 2f(x_{l},z_{l}, \rho)] \\
   & \leq \frac{1}{\mu k} \mathcal{V}(x_{0},z_{0}, \rho).
\end{align*}
Here, in the last equality, we used the identity from \eqref{als:ls:45}. This proves the first part of the theorem. Similarly, we have
\allowdisplaybreaks{
\begin{align*}
    \E[\|\bar{x}_k-x^*\|^2_{\E_{\mathcal{D}}[Z]}] & + \frac{1}{\rho} \E[\|\bar{z}_k\|_{\E_{\mathcal{D}}[W]}^2] \nonumber \\
    & \leq  \E \left[\Big \| \frac{1}{k} \sum \limits_{l=0}^{k-1} \left(x_l-x^*\right)\Big \|^2_{\E_{\mathcal{D}}[Z]}\right] +  \frac{1}{\rho}\E \left[\Big \| \frac{1}{k} \sum \limits_{l=0}^{k-1} z_l\Big \|^2_{\E_{\mathcal{D}}[W]}\right] \\
   & \leq  \E \left[\frac{1}{k} \sum \limits_{l=0}^{k-1} \| x_l-x^*  \|^2_{\E_{\mathcal{D}}[Z]}\right] + \frac{1}{\rho} \E \left[\frac{1}{k} \sum \limits_{l=0}^{k-1}  \| z_l \|^2_{\E_{\mathcal{D}}[W]}\right] \\
   & = \frac{1}{k} \sum \limits_{l=0}^{k-1} f(x_{l},z_{l}, \rho)]  \leq \frac{1}{k} \mathcal{V}(x_{0},z_{0}, \rho).
\end{align*}}
Here, in the last equality, we used the identity from \eqref{als:ls:45}. This proves the second part of Theorem \ref{als:ls:thm:3}.
\end{proof}

\section{Sampling Techniques}
\label{als:ls:sec:sampl}
In this section, we discuss several sampling strategies by which the random sketching matrix $S$ will be chosen at each iteration. Throughout this section, we will denote $S_i$ as the sketching matrix and $\E[\cdot \ | \ S_i \sim \mathcal{D} ] = \E_{i}[\ \cdot \ ]$ as the corresponding expectation.

 \begin{lemma}
\label{als:ls:lem:8}
(Lemma 2.1 in \cite{haddock:2017})  Let $\{x_k\},\ \{y_k\}$ be real non-negative sequences such that $x_{k+1} > x_k > 0$ and $y_{k+1} \geq y_k \geq 0$, then
\begin{align*}
  \sum\limits_{k=1}^{n} x_k y_k \ \geq \   \sum\limits_{k=1}^{n} \overline{x} y_k, \quad \text{where} \ \ \overline{x} = \frac{1}{n}\sum\limits_{k=1}^{n} x_k.
\end{align*}
\end{lemma}

\paragraph{Convenient Sketching.} First, we discuss a specific discrete sketching rule which is frequently used in the broader context of the SP method. This sketching rule allows us to synthesize a wide variety of sketching rules used in the literature. In the next lemma, we provide an estimation of the constant $\mu$ with respect to this specific sketching rule.

\begin{lemma}
\label{als:ls:lem:9}
Assume, $G = I$ and let $S_i$ be the sketching matrix selected with probability $p_i$,
\begin{align}
\label{als:ls:46}
   p_i = \frac{1+\rho\textbf{Tr}(S_i^\top AB^{-1}A^\top S_i)}{1+\rho\|B^{-\frac{1}{2}}A^\top \mathbf{S}\|^2_F}.
\end{align}
Denote $\mathbf{S} = [S_1,S_2,...,S_q]$. If $\mathbf{S}^\top A$ is a full rank matrix then we have
\begin{align}
\label{als:ls:47}
    \mu_c =  \frac{\rho \lambda_{\min}^+(\mathbf{S}^\top AB^{-1}A^\top \mathbf{S})}{1+\rho\|B^{-\frac{1}{2}}A^\top \mathbf{S}\|^2_F}.
\end{align}
\end{lemma}

\begin{proof}
Denote $t_i = 1+ \rho \textbf{Tr}(S_i^\top AB^{-1}A^\top S_i)$ and $t = 1+\rho\|B^{-\frac{1}{2}}A^\top \mathbf{S}\|^2_F$. Define, the following diagonal block matrices:
\begin{align*}
   & \mathbf{D}_1 = \textbf{diag}[\sqrt{t_1}H_1^{\frac{1}{2}}, \sqrt{t_2}H_2^{\frac{1}{2}},...,\sqrt{t_q}H_q^{\frac{1}{2}}], \\
   & \mathbf{D}_2 = \textbf{diag}[\sqrt{t_1}L_1^{\frac{1}{2}}, \sqrt{t_2}L_2^{\frac{1}{2}},...,\sqrt{t_q}L_q^{\frac{1}{2}}],
\end{align*}
where $L_i := \left(\frac{1}{\rho} I +B^{-\frac{1}{2}} A^\top SS^\top A B^{-\frac{1}{2}}\right)^{-1} = B^{\frac{1}{2}} F_i B^{\frac{1}{2}}.$ Using these definitions we get
\begin{align}
\label{als:ls:48}
    \E_i[B^{-\frac{1}{2}} Z B^{-\frac{1}{2}}] & = \sum\limits_{i=1}^{q} p_i B^{-\frac{1}{2}} A^\top S_i H_i S^\top _iA  B^{-\frac{1}{2}} \nonumber \\ 
    & = \frac{1}{t} B^{-\frac{1}{2}} A^\top   \sum\limits_{i=1}^{q} t_i S_i H_i S^\top _i A  B^{-\frac{1}{2}} \nonumber \\
    & = \frac{1}{t} B^{-\frac{1}{2}} A^\top  \mathbf{S} (\mathbf{D}_1)^2  \mathbf{S}^\top  A  B^{-\frac{1}{2}}.
\end{align}
In a similar fashion, we get
\begin{align}
\label{als:ls:49}
    \E_i[G^{-\frac{1}{2}} W G^{-\frac{1}{2}}] & = \sum\limits_{i=1}^{q} p_i S^\top _iA  B^{-\frac{1}{2}} B^{\frac{1}{2}} F_i B^{\frac{1}{2}} B^{-\frac{1}{2}} A^\top S_i  \nonumber \\
    & =   \sum\limits_{i=1}^{q} t_i S_i^\top A B^{-\frac{1}{2}}  L_i B^{-\frac{1}{2}}  A^\top S_i  \nonumber \\
    & = \frac{1}{t}  \mathbf{S}^\top  A  B^{-\frac{1}{2}} (\mathbf{D}_2)^2 B^{-\frac{1}{2}} A^\top  \mathbf{S}. 
\end{align}
Denote, $H_i = \left(\frac{1}{\rho} G + S_i^{\top} A B^{-1}A^\top S_i\right)^{-1}$ and $K_i = S^{\top}_i A B^{-\frac{1}{2}}$. Since, $G = I$, we can simplify matrices $H_i^{-1}$ and $L_i^{-1}$ as follows:
\begin{align*}
  & H_i^{-1} =   \frac{1}{\rho} G + S_i^{\top}A B^{-1}A^{\top} S_i = \frac{1}{\rho} I + K_i K_i^{\top} \\
  & L_i^{-1} = \frac{1}{\rho} I +B^{-\frac{1}{2}} A^\top SS^\top A B^{-\frac{1}{2}} =   \frac{1}{\rho} I + K_i^\top K_i
\end{align*}
Using the above simplified expressions, we have the following:
\begin{align*}
   \lambda_{\max}(H^{-1}_i) = \lambda_{\max}(L^{-1}_i) \leq t_i /\rho. 
\end{align*}
Now, we have
\begin{align}
\label{als:ls:50}
    \lambda_{\min}^+ (\mathbf{D}_1^2) = \lambda_{\min}^+ (\mathbf{D}_2^2) =  \min_i \left\{ \frac{t_i}{\lambda_{\max}(H_i) } \right\} \geq \rho.
\end{align}
Then, from \eqref{als:ls:48} and \eqref{als:ls:50}, we have 
\begin{align}
\label{als:ls:51}
  \lambda_{\min}^+(\E_i[B^{-\frac{1}{2}} Z B^{-\frac{1}{2}}]) & = \frac{1}{t} \lambda_{\min}^+\left(B^{-\frac{1}{2}} A^\top  \mathbf{S} (\mathbf{D}_1)^2  \mathbf{S}^\top  A  B^{-\frac{1}{2}}\right) \nonumber \\
  & = \frac{1}{t} \lambda_{\min}^+\left( \mathbf{S}^\top  A  B^{-1} A^\top  \mathbf{S} (\mathbf{D}_1)^2  \right) \nonumber \\
  & \geq \frac{1}{t} \lambda_{\min}^+\left( \mathbf{S}^\top  A  B^{-1} A^\top  \mathbf{S} \right)  \lambda_{\min}^+ (\mathbf{D}_1^2) \nonumber \\
  & \geq \frac{ \rho \lambda_{\min}^+\left( \mathbf{S}^\top  A  B^{-1} A^\top  \mathbf{S} \right)}{1+\rho\|B^{-\frac{1}{2}}A^\top \mathbf{S}\|^2_F}.
\end{align}
Similarly, considering  \eqref{als:ls:49} and \eqref{als:ls:50} we have
\begin{align}
\label{als:ls:52}
  \lambda_{\min}^+(\E_i[G^{-\frac{1}{2}} W G^{-\frac{1}{2}}]) &  = \frac{1}{t} \lambda_{\min}^+\left(\mathbf{S}^\top  A  B^{-\frac{1}{2}} (\mathbf{D}_2)^2 B^{-\frac{1}{2}} A^\top  \mathbf{S} \right) \nonumber \\
  & \geq \frac{1}{t} \lambda_{\min}^+\left( \mathbf{S}^\top  A  B^{-1} A^\top  \mathbf{S} \right)  \lambda_{\min}^+ (\mathbf{D}_2^2) \nonumber \\
  & \geq \frac{\rho \lambda_{\min}^+\left( \mathbf{S}^\top  A  B^{-1} A^\top  \mathbf{S} \right)}{1+\rho\|B^{-\frac{1}{2}}A^\top \mathbf{S}\|^2_F}.
\end{align}
Combining \eqref{als:ls:51} and \eqref{als:ls:52}, we get the required result.
\end{proof}
For the next two sketching rules, we will derive a simpler assumption than the one provided in Assumption \ref{als:ls:a:1}. First, let us define the following matrices:
\begin{align*}
   &  \mathcal{Z} = \sum\nolimits_{j=1}^{q} B^{-\frac{1}{2}}Z_{j} B^{-\frac{1}{2}}, \ \ \mathcal{W} = \sum\nolimits_{j=1}^{q} G^{-\frac{1}{2}}W_{j} G^{-\frac{1}{2}}.
\end{align*}

\begin{assumption}
\label{als:ls:a:2}
Assume sketching matrices $S_i$ satisfies the condition $\mathcal{W} \succ 0$. Moreover, assume one of the following holds:
\begin{align*}
\text{Either} \quad  \sum \limits_{i =1}^{q} S_iH_iS_i^\top  \succ \mathbf{0} \quad \text{or} \quad \textbf{Null}\left(\sum \limits_{i =1}^{q} S_iH_iS_i^\top  \right) \subseteq \textbf{Null}\left(A^\top  \right).
\end{align*}
\end{assumption}
The above assumption allows us to reformulate the original linear system of \eqref{als:ls:1} into an equivalent stochastic optimization problems like the ones provided in \eqref{als:ls:23} and \eqref{als:ls:24} under the following greedy sketching rules.

\paragraph{Greedy Sketching.} Now, we will introduce a greedy version of the above sketching rule.  Consider the following sketching rule: at iteration $k$, construct sketching set $\phi_k(\tau)$ by selecting $\tau$ sketching matrices uniformly at random from $\mathcal{S}(q) = \{S_1,S_2,...,S_q\}$, then select $S_i$ as follows:
\begin{align}
\label{als:ls:53}
i = \argmax_{i \in \phi_k(\tau)} f_i(x_{k},z_{k}, \rho).
\end{align}
Where the $i^{th}$ \textit{Sketched Loss} is defined as follows:
\begin{align}
\label{als:ls:54}
& f_i(x_{k},z_{k}, \rho) =  \frac{1}{2}\|x_k-x^*\|^2_{Z_i} + \frac{1}{2 \rho} \|z_k\|_{W_i}^2,    
\end{align}
where, the matrices $H_i, \ W_i$ can be found from Definition \ref{als:ls:d:2} with the substitution $S$ with $S_i$. The expectation calculation for this specific sketching rule is tricky. For ease of analysis, let us sort the sketched losses $f_i(x_{k},z_{k}, \rho)$ from smallest to largest. Define, $f_{\underline{\mathbf{i_j}}}(x_k,z_k, \rho) $ as the $(\tau+j)^{th}$ entry on the sorted list, i.e.,
\begin{align}
\label{als:ls:55}
  \underbrace{f_{\underline{\mathbf{i_0}}}(x_k,z_k, \rho)}_{\tau^{th}} \ \leq ... \leq \ \underbrace{f_{\underline{\mathbf{i_j}}}(x_k,z_k, \rho)}_{(\tau+j)^{th}} \ \leq ... \leq \ \underbrace{f_{\underline{\mathbf{i_{q-\tau}}}}(x_k,z_k, \rho)}_{q^{th}}.
\end{align}
Since every entry of the above list has a equal probability of getting selected, i.e., $p_i = 1/\binom{q}{\tau}$, then we can write the following:
{\allowdisplaybreaks
\begin{align}
\label{als:ls:56}
 \E_i[f_i(x_{k},z_{k}, \rho)]  & = \frac{1}{ \binom{q}{\tau}} \sum\limits_{j = 0}^{q-\tau} \binom{\tau-1+j}{\tau-1} f_{\underline{\mathbf{i_j}}}(x_{k},z_{k}, \rho).
 \end{align}}
This sketching rule can be interpreted as sketching with greedy sampling. A variant of this sketching rule has been introduced by \cite{haddock:2017} in the context of connecting two well-known sketching rules such as uniform sampling and maximum distance sampling. In their work, Morshed \textit{et. al} \cite{morshed:sketching} incorporated this sketching rule into the SP method. In the following, we will propose a lemma regarding the calculation of constant $\mu$.

\begin{lemma}
\label{als:ls:lem:10}
If we choose sketching matrix $S_i$ following the greedy sketching rule, then we have the following:
\begin{align*}
  \E_i[f_i(x_{k},  z_{k}, \rho)] \ \geq \ \frac{ \lambda_{\min}^+\left(\mathcal{Z}\right)}{2  q}  \|x_k-x^*\|^2_B + \frac{ \lambda_{\min}^+\left(\mathcal{W}\right)}{2 \rho q}  \|z_k\|^2_{G} \ \geq \  \mu_{g} \ \mathcal{V}(x_{k},z_{k}, \rho),
\end{align*}
for any iterates $x_k $ and $z_k$. Where, $\mu_g = \frac{1}{q} \min \{\lambda_{\min}^+\left(\mathcal{Z}\right), \lambda_{\min}^+\left(\mathcal{W}\right)\}$.
\end{lemma}

\begin{proof}
From the expectation definition of \eqref{als:ls:56}, we have
{\allowdisplaybreaks
\begin{align}
\label{als:ls:57}
 \E_i[f_i(x_{k},  z_{k}, \rho)] & =  \frac{1}{2 \rho \binom{q}{\tau}} \sum\limits_{j = 0}^{q-\tau} \binom{\tau-1+j}{\tau-1} \left[ \rho \|x_k-x^*\|_{Z_{\underline{\mathbf{i_j}}}}^2 + \|z_k\|_{W_{\underline{\mathbf{i_j}}}}^2  \right] \nonumber \\
 & \overset{\text{Lemma} \ \ref{als:ls:lem:8}}{\geq} \frac{1}{2 \rho \binom{q}{\tau}} \sum\limits_{j = 0}^{q-\tau}  \frac{\sum\limits_{l = 0}^{q-\tau} \binom{\tau-1+l}{\tau-1}}{q-\tau +1} \left[ \rho \|x_k-x^*\|_{Z_{\underline{\mathbf{i_j}}}}^2 + \|z_k\|_{W_{\underline{\mathbf{i_j}}}}^2  \right]  \nonumber \\
 & = \frac{1}{2 \rho (q-\tau+1)} \sum\limits_{j = 0}^{q-\tau} \left[ \rho \|x_k-x^*\|_{Z_{\underline{\mathbf{i_j}}}}^2 + \|z_k\|_{W_{\underline{\mathbf{i_j}}}}^2  \right]   \nonumber \\
 & \geq \frac{\min \left\{1, \frac{q-\tau+1}{q-s_k}\right\}}{2 \rho (q-\tau+1)}  \sum\limits_{j = 1}^{q} \left[ \rho \|x_k-x^*\|_{Z_j}^2 + \|z_k\|_{W_j}^2  \right]  \nonumber \\
  & =  \frac{1}{2 \rho q } \sum\limits_{j = 1}^{q} \left[ \rho \|x_k-x^*\|_{Z_j}^2 + \|z_k\|_{W_j}^2  \right]  \nonumber \\
 &   \geq   \frac{ \lambda_{\min}^+\left(\mathcal{Z}\right)}{2  q}  \|x_k-x^*\|^2_B + \frac{ \lambda_{\min}^+\left(\mathcal{W}\right)}{2 \rho q}  \|z_k\|^2_{G} \nonumber \\
 & \geq \ \frac{\mu_g}{2} \ \mathcal{V}(x_k, z_k, \rho),
 \end{align}}
here, we denote $\mu_g = \frac{1}{q} \min \{\lambda_{\min}^+\left(\mathcal{Z}\right), \lambda_{\min}^+\left(\mathcal{W}\right)\}$. Moreover, the last inequalities follow from Lemma \ref{als:ls:lem:3} along with Assumption \ref{als:ls:a:2}. The quantity $s_k$ is defined as  $s_k = q- \|f_k \|_0, \ f_k = [f_1(x_k,z_k, \rho),f_2(x_k,z_k, \rho),...,f_q(x_k,z_k, \rho)]^\top $, where $\|\cdot\|_0$ denotes the zero norm.
\end{proof}

\paragraph{Capped Sketching.}
The following sketching rule is also frequent in the row-action methods. Consider two sampled sketching matrix sets of sizes $\tau_1$ and $\tau_2$ respectively selected uniformly at random from $\mathcal{S}(q)$. Then, define $\mathcal{T}_k$ as follows:
\begin{align}
\label{als:ls:58}
    \begin{split}
     \mathcal{T}_k & = \left\{ i \ | \ f_i(x_k,z_k, \rho) \geq \theta  \E_i[f_i(x_k,z_k, \rho)] + (1-\theta) \E_i[f_i(x_k,z_k, \rho)] \right\}, 
    \end{split}
\end{align}
with $ \theta \in [0,1]$. Then, select $i \in \mathcal{T}_k$ with arbitrary probability. It can be easily checked that $\max_{i} f_i(x_k,z_k, \rho)  \in \mathcal{T}_k$. That implies $\mathcal{T}_k$ is not empty. The expectation with respect to the above sketching rule can be expressed as follows:
\begin{align}
\label{als:ls:59}
    \E_i[f_i(x_k,z_k, \rho) ] & = \sum \limits_{j \in \mathcal{T}_k} p_j f_j(x_k,z_k, \rho),
\end{align}
where, $0 \leq p_j \leq 1$.

\begin{lemma}
\label{als:ls:lem:11}
If we choose sketching matrix $S_i$ following the capped rule, then we have the following:
\begin{align*}
  \E_i[f_i(x_{k},  z_{k}, \rho)] \ \geq \ \frac{ \lambda_{\min}^+\left(\mathcal{Z}\right)}{2  q}  \|x_k-x^*\|^2_B + \frac{ \lambda_{\min}^+\left(\mathcal{W}\right)}{2 \rho q}  \|z_k\|^2_{G} \ \geq \  \mu_{g} \ \mathcal{V}(x_{k},z_{k}, \rho), 
\end{align*}
for any iterates $x_k $ and $z_k$. Where, $\mu_g = \frac{1}{q} \min \{\lambda_{\min}^+\left(\mathcal{Z}\right), \lambda_{\min}^+\left(\mathcal{W}\right)\}$.
\end{lemma}

\begin{proof}
First, note that from the expectation calculation, we have
 \begin{align}
 \label{als:ls:60}
   \E_i[f_i(x_k, & z_k,  \rho) ]  = \sum \limits_{j \in \mathcal{W}_k} p_j f_j(x_k,z_k, \rho) \nonumber \\
    & \geq \theta  \E[f_j(x_k,z_k, \rho) \ | \ j \sim \mathcal{G}(\tau_1)] + (1-\theta) \E[f_j(x_k,z_k, \rho) \ | \ j \sim \mathcal{G}(\tau_2)] \nonumber \\
    & \geq \ \theta  \mu_g \ \mathcal{V}(x_{k},z_{k}, \rho)  + (1-\theta)  \mu_g \ \mathcal{V}(x_{k},z_{k}, \rho)  \nonumber \\
    & = \ \mu_g \ \mathcal{V}(x_{k},z_{k}, \rho).
\end{align}
\end{proof}

\section{Special Cases}
\label{als:ls:sec:sp}
In this section, we will derive several algorithms as special cases of the proposed PS and ALS method. For ease of presentation, we will skip the special methods based on the penalty sketching method. Since, from the update formula of the ALS method, we can deduce that with the choice $z_k = 0$ in the update formula of $x_{k+1}$, we get the PS method.

\paragraph{Iterative Refinement:} Consider solving the system $Ax = b$ with $A$ positive definite matrix. Now, take $S_i= A^l, B = I$ for some $l > 0$, then form the SP update we get~\eqref{als:ls:3}:
\begin{align}
& x_{k+1} = x_k- A^{-1} (Ax_k-b) \nonumber \\
\Longleftrightarrow \  & Ay_k = b-Ax_k, \ x_{k+1} = x_k+y_k. 
\end{align}
this is precisely the well-known Wilkinson's iterative refinement method.

\paragraph{Penalty Iterative Refinement.} Take, $S_i=A^l, \ G = A^{2l+1}, \ B = I$ for some $l > 0$, then form the PS update we get:
\begin{align}
& x_{k+1} = x_k- \left(A+ \frac{1}{\rho}I\right)^{-1} (Ax_k-b) \nonumber \\
\Longleftrightarrow \ &  \left(A+ \frac{1}{\rho}I\right) y_k = b-Ax_k, \ x_{k+1} = x_k+y_k.
\end{align}
we will call this as the penalty iterative refinement method.

\paragraph{Augmented Iterative Refinement.} Take, $S_i=A^l, \ G = A^{2l+1}, \ B = I$ for some $l > 0$, then form the ALS update we get:
\begin{align*}
& z_{k+1} =  \left(A+ \frac{1}{\rho}I\right)^{-1} \left[Ax_k-b+ \frac{1}{\rho} z_k\right] , \quad   x_{k+1} = x_k- z_{k+1}. 
\end{align*}
Merging the above updates, we get the following:
\begin{align*}
& x_{k+1} = x_k- \left(A+ \frac{1}{\rho}I\right)^{-1} \left[Ax_k-b \right] + \frac{1}{\rho} \left(A+ \frac{1}{\rho}I\right)^{-1}  (x_k-x_{k-1}). 
\end{align*}
Rearranging the above with the substitution $y_k = -z_{k+1}$, we get
\begin{align}
    \left(A+ \frac{1}{\rho}I\right) y_k = b-Ax_k + \frac{1}{\rho} (x_k-x_{k-1})  , \ x_{k+1} = x_k+y_k.
\end{align}
we call this as the augmented iterative refinement method.

\paragraph{Augmented Stochastic Descent.} Take, $S_i= s_i, \ G = 1$. Then the proposed ALS method resolves into the following update formula:
\begin{align*}
& z_{k+1} =  \frac{\rho \ s_{i}^{T}(Ax_k -b)+ z_k}{1+\rho \|A^\top s_i\|^2_{B^{-1}}} , \quad   x_{k+1} = x_k- z_{k+1} B^{-1} A^\top s_i. 
\end{align*}
In that case, we have the following:
\begin{align*}
  & B^{-\frac{1}{2}} Z B^{-\frac{1}{2}}_i = B^{-\frac{1}{2}}Z_i  B^{-\frac{1}{2}} = \frac{\rho \ B^{-\frac{1}{2}} A^\top  s_i s_{i}^{T}A B^{-\frac{1}{2}}}{1+\rho \|A^\top s_i\|^2_{B^{-1}}}, \\
   & G^{-\frac{1}{2}} W_i G^{-\frac{1}{2}} = W_i = \frac{\rho \ s_i^\top  A B^{-1} A^\top  s_{i}}{1+\rho \|A^\top s_i\|^2_{B^{-1}}},
\end{align*}
and $ \mathcal{Z}  =  \sum\nolimits_{i=1}^{q} B^{-\frac{1}{2}} Z_i B^{-\frac{1}{2}} , \ \mathcal{W} =  \sum\nolimits_{i=1}^{q} G^{-\frac{1}{2}} W_i G^{-\frac{1}{2}}$. With some easy calculation, we get
\begin{align*}
    \min \{\lambda_{\min}^+(\mathcal{Z}),\lambda_{\min}^+(\mathcal{W})\} \geq \frac{\rho \lambda_{\min}^+(\mathbf{S}^\top AB^{-1}A^\top \mathbf{S})}{1+\rho\|B^{-\frac{1}{2}}A^\top \mathbf{S}\|^2_F}.
\end{align*}
Using these info we can calculate the spectral constants as follows:
\begin{align*}
& \textbf{Greedy:} \ \ i = \argmax_{j \in \phi_k(\tau)} \frac{ \|A^\top s_j\|^2_{B^{-1}} |z_k|^2  + \rho   |s_j^\top (Ax_k-b)  |^2}{  1+\rho \ \|A^\top s_i\|^2_{B^{-1}}}, \\
& \quad \quad \quad \quad \quad \quad \quad \quad \mu_g = \frac{\rho \lambda_{\min}^+(\mathbf{S}^\top AB^{-1}A^\top \mathbf{S})}{q(1+\rho\|B^{-\frac{1}{2}}A^\top \mathbf{S}\|^2_F)},   \\
  & \textbf{Convenient:} \ \ p_i = \frac{1+\rho \ \|A^\top s_i\|^2_{B^{-1}} }{1+\rho\|B^{-\frac{1}{2}}A^\top \mathbf{S}\|^2_F}, \quad   \mu_c = \frac{\rho \lambda_{\min}^+(\mathbf{S}^\top AB^{-1}A^\top \mathbf{S})}{1+\rho\|B^{-\frac{1}{2}}A^\top \mathbf{S}\|^2_F},
\end{align*}
where, $\phi_k(\tau)$ denotes the collection of $\tau$ vectors chosen uniformly at random from the sketching vector set $\mathbf{S}(q) = \{s_1, s_2,...,s_q\}$. Note that, with the choice $\rho \rightarrow \infty$ the above method with convenient sampling resolves into the \textit{Stochastic Descent} method proposed in \cite{gower:2015,richtrik2017stochastic,gower2019adaptive}. Similarly, if we take $\rho \rightarrow \infty$ in the above method with respect to the greedy sampling this method resolves into the greedy \textit{Stochastic Descent} method proposed in \cite{morshed2020stochastic}.

\paragraph{Augmented Kaczmarz} Take, $q = m, \ B =I, \ S= e_i, \ G = 1$. Then the proposed ALS method resolves into the following update formula:
\begin{align*}
z_{k+1}=  \frac{ \rho  (a_i^\top  x_k -b_i) + z_k }{1+\rho \ \|a_i\|^2}, \quad x_{k+1} =  x_k - z_{k+1} a_i.
\end{align*}
The convergence depends on the choice of the sampling process. If we select the index $i$ based on sampling rules discussed earlier we get the following complexity results:
\begin{align*}
& \textbf{Greedy:} \ \ i = \argmax_{j \in \phi_k(\tau)} \frac{ \|a_j \|^2 \ |z_k|^2  + \rho \  |a_j^\top x_k-b_j  |^2}{  1+\rho \ \|a_j \|^2}, \quad \mu_g = \frac{\rho \lambda_{\min}^+(A^\top A)}{m(1+\rho\|A\|^2_F)},  \\
  & \textbf{Convenient:} \ \ p_i = \frac{1+\rho \ \|a_{i}\|^2}{1+\rho \ \|A\|^2_F}, \quad   \mu_c = \frac{\rho \lambda_{\min}^+(A^\top A)}{1+\rho\|A\|^2_F}, 
\end{align*}
where, $\phi_k(\tau)$ denotes the collection of $\tau$ rows chosen uniformly at random out of $m$ rows of the constraint matrix $A$. The above methods with the unfiorm sampling rule have been proposed in \cite{als:kac}. Note that, the above method with the greedy sketching rule provides us the new method namely the \textit{Augmented Sampling Kaczmarz Motzkin} (ASKM) method. Furthermore, with the choice $\rho \rightarrow \infty$ the ASKM method can be sought as the \textit{Sampling Kaczmarz Motzkin} (SKM) method proposed in \cite{haddock:2019, morshed2020stochastic}. Moreover, with sample size choices $\tau =1$ and $\tau = m$ in the greedy sketching rule we recover new methods such as \textit{Augmented Randomized Kaczmarz} (ARK) method  and \textit{Augmented Motzkin} (AM) method, respectively. Finally, the first sketching rule and the second sketching rule converge into the same methods for the choices $\tau =1, \ m$ along with $\|a_i\|^2 =1, \ \forall i$.

\paragraph{Augmented Spectral Descent.} Take, $ B =A \succ 0, \ S= u_i, \ G = 1, \ q = m=n$. Then, the proposed ALS method resolves into the following update formula:
\begin{align}
    z_{k+1} & = \frac{ \rho (\xi_i u_i^\top x_k-u_i^\top b)+ z_k }{1+\rho \ \xi_i}  , \quad  x_{k+1} =  x_k - z_{k+1}    u_i,
\end{align}
here, $(\xi_i,u_i)$ is the eigenvalue-eigenvector pair for the matrix $A$.
 For this specific method, we get the following complexity results:
\begin{align*}
& \textbf{Greedy:} \ \ i = \argmax_{j \in \phi_k(\tau)} \frac{ \xi_j \ |z_k|^2  + \rho \  |\xi_j u_j^\top x_k-u_j^\top b |^2}{  1+\rho \ \xi_i}, \quad \mu_g =  \frac{\rho \lambda_{\min}^+(A)}{n(1+\rho \sum_{j} \xi_j)}, \\
  & \textbf{Convenient:} \ \ p_i = \frac{1+\rho \ \xi_i}{1+\rho \ \sum_{j} \xi_j}, \quad   \mu_c =   \frac{\rho \lambda_{\min}^+(A)}{1+\rho \sum_{j} \xi_j},
\end{align*}
where, $\phi_k(\tau)$ denotes the collection of $\tau$ eigenvectors chosen uniformly at random out of $m$ orthonormal eigenvectors of the constraint matrix $A$. Note that, the above method with the greedy sketching rule provides us a new method namely the \textit{Augmented Sampling Spectral Descent} (ASSD) method. By taking $\rho \rightarrow \infty$ in the ASSD method we get a new method namely \textit{Sampling Spectral Descent} (SSD) method. In a similar fashion, if we take $\rho \rightarrow \infty$ in the above method along with the second sketching rule, we get the so-called \textit{Spectral Descent} method proposed in \cite{stochasticspectral}.

\paragraph{Augmented Block Kaczmarz.} 
Take, $S = I_{:C}, \ B = I, \ G = I$, where $I_{:C}$ denotes the column sub-matrix of the $m \times m$ identity matrix indexed by a random set $C$.  Then the proposed ALS method resolves into the following update formula:
\begin{align*}
z_{k+1} = \left( \frac{1}{\rho} I+  A_{C:}A_{C:}^\top \right)^{-1} [A_{C:}x_k-b_C+ \frac{1}{\rho} z_k], \quad x_{k+1} = x_k-A_{C:}^\top  z_{k+1}.
\end{align*}
With the above parameter choice, we can deduce the following:
\begin{align*}
   &  B^{-\frac{1}{2}} Z B^{-\frac{1}{2}} = A_{C:}^\top  \left( \frac{1}{\rho} I+  A_{C:}A_{C:}^\top \right)^{-1} A_{C:}, \\
    & G^{-\frac{1}{2}} W G^{-\frac{1}{2}} =  A_{C:} \left( \frac{1}{\rho} I+   A_{C:}^\top  A_{C:} \right)^{-1} A_{C:}^\top .
\end{align*}
For this block variant of the augmented Kaczmarz method, we can derive the constant $\mu$ as $\mu = \min \{\lambda_{\min}^+(\E[B^{-\frac{1}{2}} Z B^{-\frac{1}{2}}]), \lambda_{\min}^+(\E[G^{-\frac{1}{2}} W G^{-\frac{1}{2}}])\}$. A better convergence result regarding the constant $\mu$ can be obtained by following the same analysis provided in \cite{NEEDELL:2014,NEEDELL:2015}. If we take $\rho \rightarrow \infty$ in the above method we get the so-called \textit{Block Kaczmarz} method \cite{NEEDELL:2014,NEEDELL:2015, gower:2015}. The proposed augmented method paves the way for designing and analyzing efficient block methods for a wide range of random blocks as well as penalty parameters.

\paragraph{Augmented Co-ordinate Newton Descent.} 
Assume, $A \succ 0$. Take, $S = I_{:C}, \ B = A,  \ G = I$, where $I_{:C}$ denotes the column sub-matrix of the $m \times m$ identity matrix indexed by a random set $C$.  Then the proposed ALS method resolves into the following update formula:
\begin{align*}
z_{k+1} = \left( \frac{1}{\rho} I+   I_{:C}^\top  A I_{:C}\right)^{-1}[ I_{:C}^\top  (A_{C:}x_k-b_C)+ \frac{1}{\rho} z_k], \quad    x_{k+1} = x_k-  I_{:C} z_{k+1}.
\end{align*}
The above iterative system is well defined as long as we choose nonempty $C$ with probability 1. Qu \textit{et. al} \cite{qu:2016} referred to this type of sampling as ``non-vacuous'' sampling. Note that for the proposed method, we can deduce the following:
\begin{align*}
   &  B^{-\frac{1}{2}} Z B^{-\frac{1}{2}} = A^{\frac{1}{2}}I_{:C} \left( \frac{1}{\rho} I+   I_{:C}^\top  A I_{:C}\right)^{-1} I_{:C}^\top  A^{\frac{1}{2}}, \\ 
    & G^{-\frac{1}{2}} W G^{-\frac{1}{2}} =  I_{:C}^\top  A^{\frac{1}{2}}\left( \frac{1}{\rho} I+  A^{\frac{1}{2}}I_{:C} I_{:C}^\top A^{\frac{1}{2}} \right)^{-1} A^{\frac{1}{2}} I_{:C}.
\end{align*}
Then the constant $\mu$ can be derived as $\mu = \min \{\lambda_{\min}^+(\E[B^{-\frac{1}{2}} Z B^{-\frac{1}{2}}]),  \lambda_{\min}^+(\E[G^{-\frac{1}{2}} W G^{-\frac{1}{2}}])\}$. A better convergence result regarding the constant $\mu$ can be obtained by following the same analysis provided in \cite{qu:2016}. Specifically, if we take $\rho \rightarrow \infty$ in the proposed augmented method, we get the so-called \textit{Co-ordinate Newton Descent} method proposed in \cite{qu:2016}.

\paragraph{Augmented CD-Positive Definite (ARCD-PD)} Assume, $A \succ 0$. Then, with the choice $q = m =n, \ B = A, \ S= e_i, \ G = 1$, the proposed ALS method resolves into the following: 
\begin{align*}
z_{k+1}=  \frac{ \rho  (a_i^\top  x_k -b_i) + z_k }{1+\rho \ A_{ii}}, \quad x_{k+1} =  x_k - z_{k+1} e_i,
\end{align*}
where, $A_{ii}$ is the $i^{th}$ diagonal entry of matrix $A \in \R^{n \times n}$. For this variant of RCD we get the following complexity results:
\allowdisplaybreaks{\begin{align*}
& \textbf{Greedy:} \ \ i = \argmax_{j \in \phi_k(\tau)}  \frac{ A_{jj} \ |z_k|^2  + \rho \  |a_j^\top x_k-b_j |^2}{  1+\rho \ A_{jj}}, \quad \mu_g = \frac{\rho \lambda_{\min}^+(A)}{n(1+\rho \textbf{Tr}(A))}  \\
  & \textbf{Convenient:} \ \ p_i = \frac{1+\rho \ A_{ii}}{1+\rho \ \textbf{Tr}(A)}, \quad   \mu_s = \frac{\rho \lambda_{\min}^+(A)}{1+\rho \textbf{Tr}(A)},
\end{align*}}
where, $\phi_k(\tau)$ denotes the collection of $\tau$ rows chosen uniformly at random out of $m$ rows of the constraint matrix $A$. Note that, the above method with the greedy sketching rule provides us the new method namely the \textit{Augmented Sampling CD-PD} (ASCD-PD) method. Furthermore, with the choice $\rho \rightarrow \infty$ the ASCD-PD method can be sought as the \textit{Sampling CD-PD} (SCD-PD) method proposed in \cite{morshed2020stochastic}. Moreover, with sample size choices $\tau =1$ and $\tau = m$ in the greedy sketching rule we recover new methods such as \textit{Augmented Randomized CD-PD} (ARCD-PD) method  and \textit{Augmented Motzkin CD-PD} (AMCD-PD) method, respectively. Similarly, if we take $\rho \rightarrow \infty$ in the ARCD-PD method then we get the CD method proposed in \cite{lewis:2010,gower:2015}. Finally, the first sketching rule and the second sketching rule converge into the same methods for the choices $\tau =1, \ m$ along with $A_{ii} =1, \ \forall i$.

\paragraph{Augmented RCD-Least Square (ARCD-LS)} Now, with the choice $q = n, \ B = A^\top A, \ S= Ae_i = A_i, \ G = 1$ the proposed ALS method resolves into the following:
\begin{align*}
z_{k+1}=  \frac{ \rho  A_i^\top (A x_k -b) + z_k }{1+\rho \ \|A_{i}\|^2}, \quad x_{k+1} =  x_k - z_{k+1} e_i,
\end{align*}
where, $A_i$ denotes the $i^{th}$ column of matrix $A$. Using the proposed lemmas, we get the following complexity results for the above method:
\begin{align*}
& \textbf{Greedy:} \ \ i = \argmax_{j \in \phi_k(\tau)}  \frac{ \|A_j \|^2 \ |z_k|^2  + \rho \ \big |A_j^\top (Ax_k-b) \big |^2}{  1+\rho \ \|a_j \|^2}, \quad \mu_g = \frac{\rho \lambda_{\min}^+(A^\top A)}{n(1+\rho \|A\|^2_F)}  \\
  & \textbf{Convenient:} \ \ p_i = \frac{1+\rho \ \|A_{i}\|^2}{1+\rho \ \|A\|^2_F}, \quad   \mu_s = \frac{\rho \lambda_{\min}^+(A^\top A)}{1+\rho \|A\|^2_F},
\end{align*}
where, $\phi_k(\tau)$ denotes the collection of $\tau$ rows chosen uniformly at random out of $n$ columns of the constraint matrix $A$. Note that, the above method with the greedy sketching rule provides us the new method namely the \textit{Augmented Sampling CD-LS} (ASCD-LS) method. Furthermore, with the choice $\rho \rightarrow \infty$ the ASCD-LS method can be sought as the \textit{Sampling CD-LS} (SCD-LS) method proposed in \cite{morshed2020stochastic}. Moreover, with sample size choices $\tau =1$ and $\tau = m$ in the greedy sketching rule, we recover new methods such as \textit{Augmented Randomized CD-LS} (ARCD-LS) method  and \textit{Augmented Motzkin CD-LS} (AMCD-LS) method, respectively. Similarly, if we take $\rho \rightarrow \infty$ in the ARCD-LS method then we get the CD method proposed in \cite{lewis:2010,gower:2015}.

\paragraph{Gaussian Sketches} Now, we will discuss some of the methods that can be obtained by choosing sketch vectors as Gaussian vectors. Assume sketch vector $S$ is a Gaussian vector with zero mean and covariance matrix $\Sigma $, i.e, $\zeta \sim \mathcal{N}(0,\Sigma), \ G = 1$. Then the proposed ALS method resolves into the following update formula:
\begin{align*}
& z_{k+1} =  \frac{\rho \ \zeta^{T}(Ax_k -b)+ z_k}{1+\rho \|A^\top \zeta\|^2_{B^{-1}}} , \quad   x_{k+1} = x_k- z_{k+1} B^{-1} A^\top \zeta.
\end{align*}
In that case, we have the following:
\begin{align*}
   & B^{-\frac{1}{2}} Z B^{-\frac{1}{2}} = B^{-\frac{1}{2}}Z  B^{-\frac{1}{2}} = \frac{\rho \ B^{-\frac{1}{2}} A^\top  \zeta \zeta^{T}A B^{-\frac{1}{2}}}{1+\rho \|A^\top \zeta\|^2_{B^{-1}}}, \\
   & G^{-\frac{1}{2}} W G^{-\frac{1}{2}} = W = \frac{\rho \ \zeta^\top  A B^{-1} A^\top  \zeta}{1+\rho \|A^\top \zeta\|^2_{B^{-1}}}.
\end{align*}

\paragraph{Augmented Gaussian Kaczmarz.} Take, $ B =I, \ S= \eta \sim \mathcal{N}(0,I)$. Then the proposed ALS method resolves into the following update formula:
\begin{align*}
 z_{k+1} = \frac{\rho \eta^{T}(Ax_k -b)+ z_k}{1+\rho \|A^\top \eta\|^2}, \quad x_{k+1} = x_k- z_{k+1} A^\top \eta. 
\end{align*}

\paragraph{Augmented Gaussian Least-Square.} Take, $ B = A^\top A, \ S  \sim \mathcal{N}(0,AA^\top )$. For simplification, denote $S = A \eta$ where, $\eta  \sim \mathcal{N}(0,I)$, then we get the following update:
\begin{align*}
 z_{k+1} = \frac{\rho \ \eta^{T} A^\top (Ax_k -b)+ z_k}{1+\rho \|A\eta\|^2}, \quad x_{k+1} = x_k- z_{k+1} \eta. 
\end{align*}

\paragraph{Augmented Gaussian Positive Definite.} Take, $ B =A, \ S= \eta \sim \mathcal{N}(0,I)$. Then the proposed ALS method resolves into the following update formula:
\begin{align*}
 z_{k+1} = \frac{\rho \eta^{T}(Ax_k -b)+ z_k}{1+\rho \|\eta\|^2_A}, \quad x_{k+1} = x_k- z_{k+1}  \eta. 
\end{align*}


\section{Conclusion}
\label{sec:concl}

In this work, we proposed two novel variants of the SP method for solving linear system of equations. The proposed penalty sketch and augmented Lagrangian sketch opens new avenues for the development of efficient randomized algorithms for solving LS problems. We proved global linear convergence results for the proposed algorithms that generalize existing convergence results regarding the SP method. The proposed
sketching methods outperform the existing SP methods on various synthetic and real-life sparse datasets.

\paragraph{Extension} The proposed PS and ALS approaches can be extended to several sketching methods. Following the PS and ALS approach, one can develop penalty and augmented variants of the subspace Newton method proposed in \cite{gower2019rsn}. Another obvious extension of the current work would be to design penalty/augmented randomized projection methods for the linear/convex feasibility problem \cite{necoara:2019,morshed:sketching}. Using the PS and ALS approaches proposed in this work, one can design penalty and augmented variants of the Newton-Raphson sketch that was proposed in \cite{yuan2020sketched}. Finally, following the proposed sketching approaches, one can build penalty/augmented versions of Quasi-Newton sketching proposed in \cite{gower:2017,NIPS:2018}.

\appendix

\section*{Appendix}

\begin{lemma}
\label{als:ls:lem:12}
(Lemma 14 in \cite{gower2019adaptive}) If the relation $\textbf{Null}(P) \subset \textbf{Null}(M^\top )$ holds for some matrix $M$ and positive semi-definite matrix $P \succeq 0$. Then, we have
\begin{align*}
\textbf{Range}(M^\top ) = \textbf{Range}(M^\top PM).
\end{align*}
\end{lemma}

\begin{lemma}
\label{als:ls:lem:120}
(Woodbury Matrix Identity) Assume matrices $A \in \R^{q \times q}$, and $C \in \R^{r \times r }$ are invertible. Then the following identity holds:
\begin{align}
\label{als:ls:61}
\left(A+ UCV\right)^{-1} = A^{-1} - A^{-1}  U \left(C^{-1}+VA^{-1} U\right)^{-1} VA^{-1},
\end{align}
for any matrices $U \in \R^{q \times r}$ and $V \in \R^{r \times q}$ such that $A+ UCV$ and $C^{-1}+VA^{-1} U$ are non-singular.
\end{lemma}

\begin{lemma}
\label{als:ls:lem:13}
The following relations hold:
\begin{enumerate}
    \item $ F = \rho B^{-1} - \rho B^{-1} Z B^{-1}$.
    \item $ GH G = \rho G - \rho W$.
    \item $HS^\top A B^{-1}A^\top S H = H- \frac{1}{\rho} H GH$.
    \item $A^\top S H- ZB^{-1}A^\top S H = \frac{1}{\rho} A^\top S HG H$.
    \item $ZB^{-1}Z = Z-\frac{1}{\rho} A^\top S HG HS^\top A$.
\end{enumerate}
\end{lemma}

\begin{proof}
Since, $B, G \succ 0$, by Woodburry matrix identity we have the following:
\begin{align*}
    \left(B+ \rho A^\top S G^{-1} S^\top A \right)^{-1} & = B^{-1} -  B^{-1}  A^\top  S \left( \frac{G}{\rho}+ S^\top A B^{-1} A^\top S\right)^{-1} S^\top  A B^{-1} \\
  & = B^{-1} -  B^{-1} Z B^{-1}.
\end{align*}
The first identity follows from the definition of $F$. The second identity follows by applying the Woodburry identity to the matrix $H$. Now, we have
\begin{align}
   HS^\top A B^{-1}A^\top S H & = H S^\top A B^{-1} A^\top S H =  H [H^{-1}-\frac{1}{\rho} G]H  = H- \frac{1}{\rho} H GH.
\end{align}
This proves the third identity. Now the fourth and fifth identity follows from the definition, i.e., $Z = A^\top S HS^\top A$.
\end{proof}

\section*{Acknowledgements}
The Author would like thank Rober M. Gower, Michał Dereziński, and Michael Mahoney for their thoughtful comments on an earlier version of the manuscript.


\bibliographystyle{unsrt}
\bibliography{references}


\end{document}